\documentclass[3p]{elsarticle}
\usepackage{amsmath,amsfonts,enumerate,amsthm,color, amssymb,multirow}
 \usepackage{url}
\usepackage{pgf}
\newtheorem{lemma}{Lemma}[section]
\newtheorem{proposition}[lemma]{Proposition}
\newtheorem{corollary}[lemma]{Corollary}
\newtheorem{theorem}[lemma]{Theorem}

\newtheorem{example}[lemma]{Example}

\newcommand{\eps}{\epsilon}
\newcommand{\ZZ}{\mathbb{Z}}

\DeclareMathOperator{\colors}{{\cal K}}

\newcommand{\comment}[1]{}

\DeclareMathOperator{\Cay}{Cay}
\DeclareMathOperator{\Aut}{Aut}

\DeclareMathOperator{\Prism}{Pr}
\DeclareMathOperator{\Ml}{Ml}
\DeclareMathOperator{\GP}{GP}
\DeclareMathOperator{\GPr}{GPr}

\DeclareMathOperator{\HTG}{HTG}
\DeclareMathOperator{\girth}{g}
\DeclareMathOperator{\signat}{sgn}

\title{{\bf On cubic rainbow domination regular graphs\tnoteref{t1}}}
\tnotetext[t1]{This work is supported in part by the Slovenian Research Agency (ARRS), research program P1-0285 and research projects J1-9108, J1-1694, J1-1695.}

\author{Bo\v stjan Kuzman
}
\ead{bostjan.kuzman@pef.uni-lj.si}
\address{University of Ljubljana, Faculty of Education, Dept. of Mathematics and Computer Science
}

\begin{document}

\begin{abstract}
A $d$-regular graph $X$ is called $d$-rainbow domination regular or $d$-RDR, if its $d$-rainbow domination number $\gamma_{rd}(X)$ attains the lower bound $n/2$ for $d$-regular graphs, where $n$ is the number of vertices. In the paper, two combinatorial constructions to construct new $d$-RDR graphs from existing ones are described and two general criteria for a vertex-transitive $d$-regular graph to be $d$-RDR are proven. A list of vertex-transitive 3-RDR graphs of small orders is produced and their partial classification into families of generalized Petersen graphs, honeycomb-toroidal graphs and a specific family of Cayley graphs is given by investigating the girth and local cycle structure of these graphs.
\end{abstract}


\maketitle


\section{Introduction}
\label{sec:intro}

The concept of rainbow domination extends the classical domination problem in graphs and was initially introduced by Hartnel and Rall \cite{Hartnell Rall 2004} in the context of the famous Vizing conjecture on domination number of cartesian products. We define the rainbow domination number of a graph as follows. Let $G=(V,E)$ be a finite, simple graph, where $V,E$ are the sets of vertices and edges, respectively. We denote adjacent vertices by $u\sim v$ and by  $N(v)=\{u\in V\colon u\sim v\}$ we denote the \emph{open neighbourhood} of the $v$.
For a given nonnegative integer $k$ we denote by $\colors=\{1,\ldots,k\}$ the set of \emph{colors}. A function $f\colon V\to 2^{\colors}$ that assigns to each vertex $v\in V$ a subset of colors $f(v)\subseteq {\colors}$, is  called 
a \emph{$k$-rainbow dominating function} (or $k$-RDF) on $G$, if
$f(v)=\emptyset$ implies $\cup_{u\in N(v)}f(u)={\colors}$ for all $v\in V$.
In other words, for every \emph{non-colored vertex} $v$ (that is, a vertex with $f(v)=\emptyset$) all possible colors appear in its neighbourhood.
For a given $k$-RDF, we then define the \emph{weight} of $f$ as 
$w(f)=\sum_{v\in V}|f(v)|.$ As defined in~\cite{Bresar H R 2008}, the \emph{$k$-rainbow domination number of $G$} is defined as
the minimum weight of all $k$-rainbow domination functions on $G$:
$$\gamma_{rk}(G)=\min\{w(f)\mid f\colon V\to2^{\colors}\text{ is a $k$-RDF}\},$$
and any $k$-rainbow dominating function $f$ on $G$ of minimal weight $w(f)=\gamma_{rk}(G)$ is called a \emph{$\gamma_{rk}(G)$-function}. 

Rainbow domination has natural applications in efficient network construction.
For a given graph $G$ and positive integer $k$, determination of the exact value $\gamma_{rk}(G)$ is known to be NP-complete~\cite{Bresar Kraner 2007, Chang 2009}, hence many research results were focused on determining some upper or lower bounds for $\gamma_{rk}$, see for instance \cite{Wang 2013}, \cite{Shao 2014}, \cite{Fujita 2015}, \cite{Furuya 2018}, or exact values for specific graph families such as trees, products of paths and cycles, generalized Petersen graphs, grid graphs, etc., see for instance \cite{Chang 2009}, \cite{Shao 2019}, \cite{Stepien 2014}, \cite{Zerovnik ckk}, or a survey paper by Bre\v sar \cite{Bresar survey}.
\medskip

In \cite{kuzman regular}, regular graphs were studied. A general lower bound for $k$-rainbow domination number of $d$-regular graphs was determined as $\gamma_{rk}(G)\geq \frac{kn}{2d}$, where $n=|V(G)|$. For $d=k$, regular graphs attaining the lower bound $\gamma_{rd}(G)=\frac{n}{2}$ are called  \emph{$d$-rainbow domination regular graphs}, or \emph{$d$-RDR graphs}.  
Specific properties and several families of $2-$, $3-$ and $4-$RDR's were identified, such as cycles, prisms, M\"obius ladders, and wreath graphs. As such graphs represent a very efficient regular network structure, in which each main node is dominated by $d$ side nodes of different types, we pursue their further study in this paper.

Although initial examples of $d$-RDR graps have suggested that they might neccesarily have specific symmetry properties such as vertex-transitivity, several examples of non-vertex-transitive $3$-RDR and $4$-RDR graphs were later described by \v Zerovnik in \cite{Zerovnik not VT}. Moreover, our own computer investigations by Magma~\cite{Magma} have shown that among all bipartite cubic graphs of small order, only a tiny share of $3$-RDR's are vertex transitive (see Table~\ref{tab: bicubic} for graphs of orders up to $96$), making the complete classification of all $3$-RDR's currently out of reach. 

However, the number of vertex-transitive bicubic graphs is much smaller, and for small orders, roughly half of them are also $3$-RDR. With many results on symmetric cubic graphs readily available, such as the census of vertex-transitive cubic graphs by Poto\v cnik, Spiga and Verret \cite{Potocnik Spiga Verret CVT} and classification of vertex-transitive cubic graphs of small girths  in \cite{Sparl Eiben Jajcay} and \cite{Potocnik Vidali girth six}, classification of all $3$-RDR's with a specific symmetry property such as vertex-transitivity is a natural project we pursue in this paper.

\begin{table}[h!]
\centering
\begin{tabular}{|c|c|c|c|c|c|}
\hline 
\multirow{2}{*}{Order} & \multirow{2}{*}{BC} & \multirow{2}{*}{3-RDR} & 
			\multirow{2}{*}{VTBC} & \multirow{2}{*}{VT 3-RDR}  & Unidentified\\ 
 						& &  & &  & {\small (girth 4/6/8/10)}\\ \hline
6 &1&1&1&1&-\\
12 & 5& 3& 2& 1&-  \\
18 & 149 &  37& 3& 2&-\\
24 & 29579& 1998& 7& 3&-\\
30 & 
       $2.3\cdot 10^7$ & .& 6& 2&- \\
36 & 
       $3.9\cdot 10^{10}$&. & 9& 5&$1/0/0/0$\\
42 & .  & . & 7 &3 &- \\
48 & . & . & 25 & 12&$3/2/0/0$\\
54 & . & .  & 9 & 5&$0/1/0/0$\\
60 & . & . & 20 & 5&$1/0/0/0$\\
66& . & .& 9 & $\geq 3$ &-\\
72& . & . & 29 &$\geq 12$& $1/3/2/0$\\ 
78& . & . & 11 & $\geq 4$&-\\
84& . & . & 24 & $\geq 9$&$1/2/2/0$ \\
90& . &. & 16& $\geq 6$ & -\\
96&.  & . & 79& $\geq 43$&$9/5/16/1$\\ \hline
\end{tabular}
\caption{Total numbers of nonisomorphic  bicubic graphs (BC), 3-rainbow domination regular graphs (3-RDR),  vertex-transitive bicubic graphs (VTBC), vertex-transitive 3-RDR graphs (VT 3-RDR) of orders $\leq 96$, with placeholders for undetermined values. Column Unidentified gives the total numbers of VT 3-RDR graphs of girths 4/6/8/10 that are not isomorphic to any $\GP$, $\HTG$ or $X_n$ graph as defined in Section~\ref{sec VT 3RDR}. Data for column BC was compiled from \cite{Coolsaet}, for column VTBC from \cite{Potocnik Spiga Verret CVT},  and data in other columns was obtained by our own computations in the Magma system~\cite{Magma}, using criteria from Theorems \ref{krit1} and \ref{krit2} for last six rows of column VT 3-RDR.}\label{tab: bicubic}
\end{table}

\begin{table}
\centering
\begin{tabular}{|c|c|c |l|}
\hline
Order & Index & Girth & Graph descriptions\\ \hline 

6  & [6,2] &4&$\Ml(3)\cong \HTG(1,6,3)\cong K_{3,3}$ \\ \hline

12&[12,3] &4& $\Prism(6)\cong\HTG(1,12,3)\cong\HTG(2,6,0)$ \\ \hline

18&[18,1] &4& $\Ml(9)\cong \HTG(1,18,3)\cong\HTG(1,18,9)$ \\
   &[18,4] &6& $\HTG(3,6,3)$ aka Pappus graph \\ \hline

24&[24,2] &6& $\HTG(2,12,6)\cong\GP(12,5)$ aka Nauru graph \\
   &[24,4] &4& $\Prism(12)\cong\HTG(1,24,3)\cong\HTG(2,12,0)$ \\
   &[24,6] &6& $\HTG(1,24,9)\cong\HTG(4,6,0)$ \\ \hline 

30 &[30,3] & 6& $\HTG(1,30,9)\cong\HTG(5,6,3)$ \\
    &[30,5] & 4& $\Ml(15)\cong \HTG(1,30,3)\cong\HTG(1,30,15)$\\ \hline 
    
36 & [36,6] & 6& $\HTG(1,36,9)\cong\HTG(1,36,15)\cong\HTG(2,18,6)$\\
    &[36,8]  & 4&$\Prism(18)\cong \HTG(1,36,3)\cong\HTG(2,18,0)$ \\
    &[36,10]& 6&$\HTG(3,12,3)\cong\HTG(6,6,0)$ \\
   &[36,11]& 6& $X_3$  \\
   &[36,12]& 4& A generalized truncation of some quartic AT graph\\ \hline
\end{tabular}
\caption{Non-isomorphic vertex-transitive 3-RDR graphs of orders $\leq 36$, where $[n,k]$ is the index of the graph
in the Census of vertex-transitive graphs~\cite{Potocnik Spiga Verret CVT}.}\label{tab: ord<36}
\end{table}

The structure of the paper is the following.
In Section 2, we collect the basic results on $d$-RDR graphs from \cite{kuzman regular} (Theorem \ref{thm: gamma rk pred}) and describe 2 combinatorial ways to construct new $d$-RDR graphs from existing ones (Propositions \ref{switch} and \ref{stitch}). In Section 3, previously known examples of infinite families of vertex-transitive $d$-RDR graphs are presented (Proposition \ref{basic RDR}) and two new group-theoretical criteria for a vertex-transitive graph to be $d$-RDR are proven (Theorems \ref{krit1} and \ref{krit2}), enabling us to efficiently compile a (possibly incomplete) list of graphs of small orders (Tables \ref{tab: bicubic} and \ref{tab: ord<36}). In Section 4, two important families of cubic graphs are then investigated in full detail. The criteria for generalized Petersen graphs $\GP(n,k)$ and honeycomb toroidal graphs $\HTG(m,n,\ell)$ to be vertex-transitive $3$-RDR graphs are described in terms of their parameters (Theorems \ref{krit GP} and \ref{krit HTG}). Vertex-transitive $3$-RDR graphs from these two classes are classified further by their girth and signature, that is, the local structure of the girth cycles (Theorems \ref{3RDR GP} and \ref{3RDR HTG sgn}). Moreover, another infinite family of vertex-transitive $3$-RDR graphs is constructed as Cayley graphs over a direct product of symmetric group $S_3$ and dihedral group $D_n$ and their overlap with $\GP$ and $\HTG$ graphs is determined (Theorem \ref{thm Xn} and Corollary \ref{Xn ni HTG}). However, it turns out that these three families do not classify the $3$-RDR graphs completely. Hence, Section 5 concludes the paper by a short discussion of several open questions.

\section{Some general observations on $d$-RDR graphs and their colorings}
First, we collect some basic facts on rainbow domination in regular graphs. 
The following Theorem combines the main results from \cite[Theorems 1.2, 3.1 and Corollary 3.5]{kuzman regular}.

\begin{figure}[h!]
\centering
\includegraphics[width=15cm]{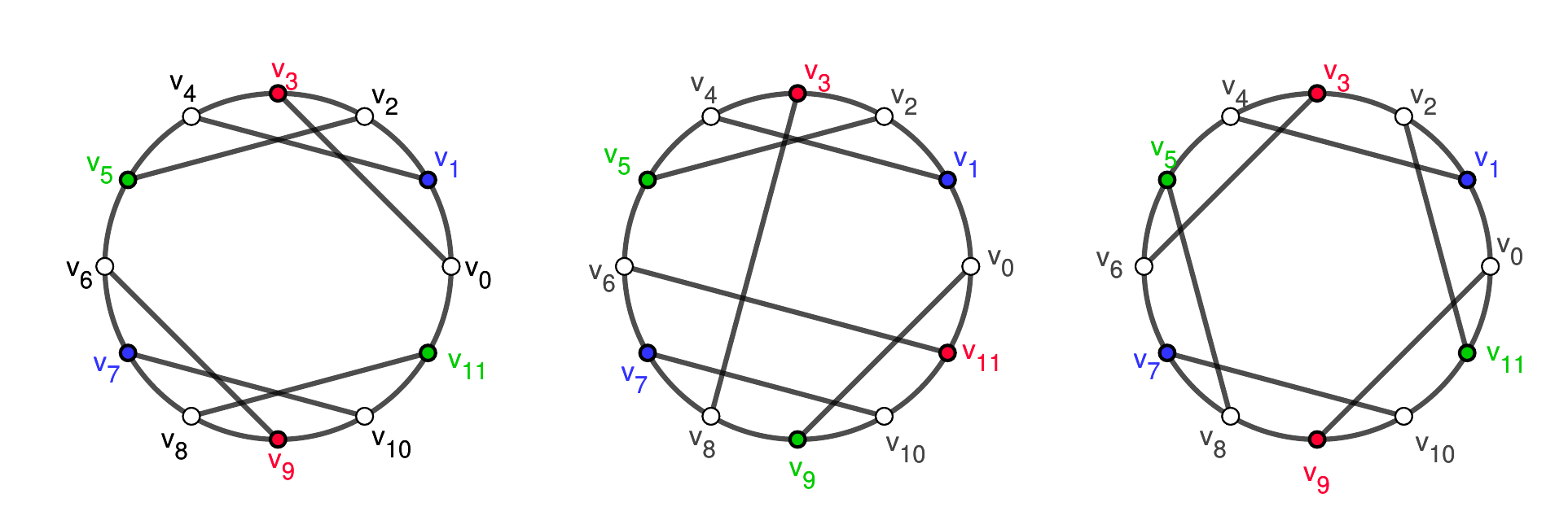}
\caption{The three non-isomorphic 3-RDR graphs of order 12. The graph on the right is vertex-transitive and isomorphic to the prism $\Prism(6)$}\label{3RDR order 12}
\end{figure}

\begin{theorem}[Kuzman, \cite{kuzman regular}]\label{thm: gamma rk pred}
Let $X$ be a $d$-regular graph of order $n$ and $k\leq 2d$. Then:
\begin{enumerate}[$(i)$]
\item 
$\gamma_{rk}(X)\geq \left\lceil \frac{kn}{2d}\right\rceil$.
\item
If $\gamma_{rk}(X)=\frac{kn}{2d}$ for some $k<2d$, then $k\geq d$, $2d|n$ and $G$ is bipartite with
bipartition sets of size $n/2$.
Moreover, for any $\gamma_{rk}(X)$-function $f$ on $X$, the sets of colored and non-colored vertices define the bipartition and each color $i\in\colors$ appears exactly $n/2d$ times (implying that each colored vertex is colored with a single color).

\item $G$ is $d$-rainbow domination regular if and only if $\gamma_{rk}(X)=\frac{kn}{2d}$ for all $d\leq k\leq 2d$.
\end{enumerate}
\end{theorem}
It follows that every $d$-RDR graph is bipartite and one of the bipartition sets can be further partitioned into $d$ subsets of same size $n/2d$, each representing a set of vertices of a single color. Every colored vertex has $d$ non-colored neighbours and every non-colored vertex has $d$ neighbours of different colors (each colored with a single color). In fact, we can view $d$-RDR property as a special case of a standard vertex-coloring with $(d+1)$ colors, where one of the colors, for instance white, has neighbours of all other colors (see Figure \ref{3RDR order 12}). 
For cubic graphs, the following necessary condition for a graph to be $3$-RDR is immediately obtained from Theorem \ref{thm: gamma rk pred}.
\begin{corollary}
Suppose that a cubic graph $X$ of order $n$ is $3$-RDR. Then $6|n$ and $X$ is bipartite with bipartition sets of size $n/2$. Moreover, for any $\gamma_{r3}(X)$-function, each color $i\in\{1,2,3\}$ appears exactly $n/6$ times, the sets of colored and non-colored vertices define the bipartition, and every colored vertex is colored with a single color.
\end{corollary}

We also point out that a $d$-RDR coloring function on a $d$-RDR graph is not necessarily unique (despite allowing permutations of colors). Moreover, the bipartition set of colored vertices is sometimes uniquely determined, but sometimes either of bipartition sets can be colored.
\begin{example}
Observe the three non-isomorphic 3-RDR graphs of order 12 in Figure~\ref{3RDR order 12}. The first one (on the left) has two essentialy different coloring functions (colors of $v_7$ and $v_{9}$ can be exchanged while preserving other colors), and both bipartition sets can be colored equivalently because of symmetry. The second graph (in the middle) has a unique 3-RDR coloring function (up to permutation of colors) and it is not difficult to check that its set of non-colored vertices is uniquely determined. The third graph (on the right) has a unique coloring function (up to permutation of colors and graph symmetry) and each of its bipartition sets can be colored or non-noncolored.
Note also that while the last graph is isomorphic to a prism and is obviously vertex-transitive, vertex-transitivity is generally not enough to obtain a unique 3-RDR coloring. As checked by computer, graph $X_3$ of order 36 (see Table~\ref{tab: ord<36} and subsection~\ref{subs Xn} for definition) is the smallest example of a vertex-transitive 3-RDR graph with two essentially different 3-RDR colorings.
\end{example}
 
\begin{figure}[h!]
\centering
\includegraphics[width=14cm]{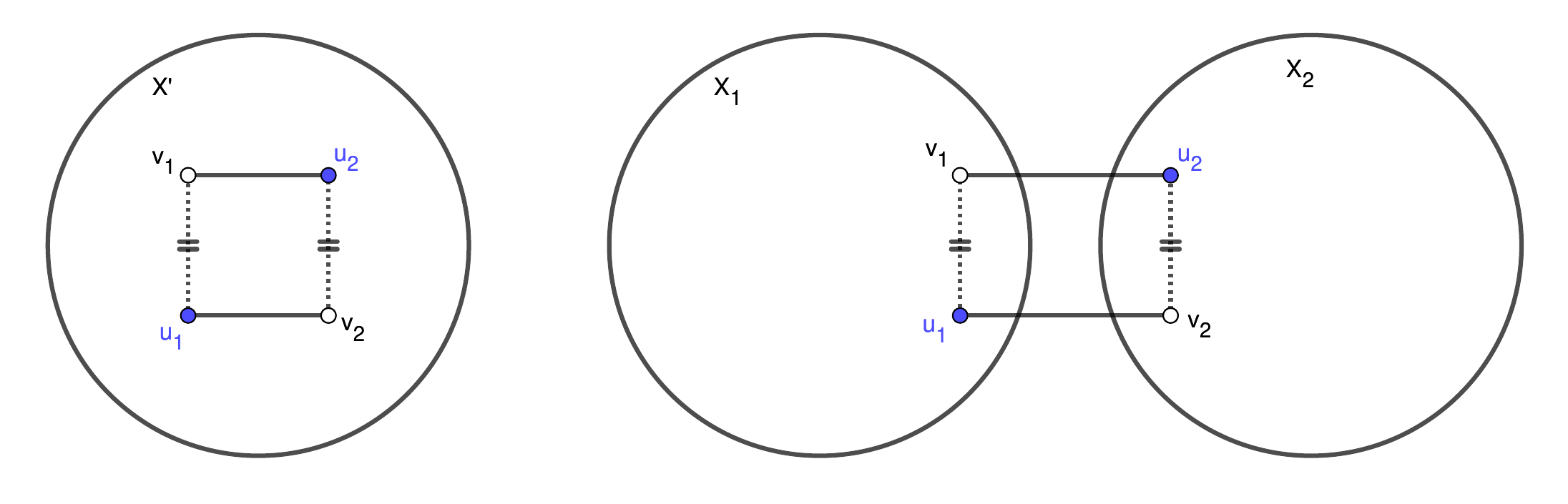}
\caption{Edge switching (left) and graph stitching (right) operations  on $d$-RDR graphs}\label{sw-st}
\end{figure}

Next, we describe two combinatorial constructions to obtain new $d$-RDR graphs from existing ones. The first one transforms a $d$-RDR graph into another one of the same order (but possibly isomorphic to the original one), while the second one 
combines two $d$-RDR graphs into a larger one (see Figure~\ref{sw-st} for the illustration of both operations). 
\begin{proposition}[RDR edge switching]\label{switch}
Let $X$ be a $d$-RDR graph and let $f$ be a $d$-RDR coloring function on $X$. Suppose that $e_1=\{u_1,v_1\}$ and $e_2=\{u_2,v_2\}\in E(X)$ are two edges with endpoints of same color: $f(u_1)=f(u_2)\ne\emptyset$ (implying $f(v_1)=f(v_2)=\emptyset$). Construct graph $X’$ with same vertex set and edge set $E(X')=(E(X)\cup\{e'_1,e'_2\})\setminus\{e_1,e_2\}$, where $e'_1=\{u_1,v_2\}$ and $e'_2=\{u_2,v_1\}$. Then $X’$ is a $d$-RDR graph with a $d$-RDR coloring function $f'=f$.
\end{proposition}
\begin{proof}
Obviously, the coloring function $f'$ is $d$-RDR on $X'$, as vertices $v_1,v_2$ have neigbours of all colors, and the rest is unchanged.
\end{proof}
\begin{example}
It is not difficult to check manually that all three non-isomorphic $3$-RDR graphs of order $12$ can be obtained from one another by edge switching operations. For instance, by switching edges $\{v_1,v_4\}$ and $\{v_7,v_{10}\}$ in the graph on the left we obtain a graph isomorphic to the graph on the right. Similarily, by switching edges $\{v_3,v_8\}$ and $\{v_6,v_{11}\}$ of the graph in the middle we obtain a graph isomorphic to the one on the left. Using computer we also checked that any pair of 3-RDR graphs of orders $18$ or $24$ can be transformed from one graph into another by a sequence of edge switching operations. One might ask whether this is true in general.
\end{example}

\begin{proposition}[RDR graph stitching]\label{stitch}
Let $X_1,X_2$ be two $d$-RDR graphs and let $f_1,f_2$ be their $d$-RDR coloring functions. Suppose that
$e_1=\{u_1,v_1\}\in E(X_1)$ and $e_2=\{u_2,v_2\}\in E(X_2)$ are two edges with endpoints of same color: $f_1(u_1)=f_2(u_2)\ne\emptyset$ (implying $f(v_1)=f(v_2)=\emptyset$).
Let $X_1*X_2$ be a graph with vertex set $V(X_1*X_2)=V(X_1)\cup V(X_2)$ and edge set $E(X_1*X_2)=(E(X_1)\cup E(X_2)\cup\{e'_1,e'_2\})\setminus
\{e_1,e_2\}$, where $e'_1=\{u_1,v_2\}$ and $e'_2=\{u_2,v_1\}$. Then $X_1*X_2$ is a $d$-RDR graph with $d$-RDR coloring function $f=f_1\cup f_2$. Moreover, similar operation can be defined with several pairs of edges with endpoints of same color in each pair.
\end{proposition}
\begin{example}
It is not difficult to check that all three non-isomorphic $3$-RDR graphs of order 12 in Figure~\ref{3RDR order 12} can be obtained by stitching 2 copies of the complete bipartite graph $K_{3,3}$ along appropriate pairs of edges (for instance, replace edges $\{v_5,v_6\}$ and $\{v_0,v_{11}\}$ in the graph on the left with edges $\{v_5,v_0\}$ and $\{v_6,v_{11}\}$).
\end{example}
\section{Vertex-transitive $d$-RDR graphs}

Graph $X$ is called \emph{vertex-transitive}, if its automorphism 
group $\Aut(X)$ acts transitively on the vertex set $V(X)$.
Recall that the \emph{prism graph} $\Prism(n)$ of order $2n$ is a cubic graph with vertex set $\{u_i,v_i\colon i\in \ZZ_n\}$ and edge set $\{\{u_i,u_{i+1}\}, \{v_i,v_{i+1}\},\{u_i,v_i\}\colon i\in\ZZ_n\}$, while a \emph{M\"obius ladder graph} is obtained from the prism by replacing two edges $\{u_{i-1},u_0\},\{v_{i-1},v_{0}\}$ with $\{u_{i-1},v_0\},\{v_{i-1},u_0\}$. \emph{Wreath graph} $W(n)$ is a quartic graph with the same vertex set and edge set $\{\{u_i,u_{i+1}\},\{u_i,v_{i+1}\},\{v_i,v_{i+1}\},\{v_i,u_{i+1}\}\colon i\in\ZZ_n\}$; it can also be described as the lexicographic product of an $n$-cycle $C_n$ and two vertices $2K_1$. Note that all graphs in these families are vertex-transitive, and the following Proposition collects several infinite families of connected vertex transitive $d$-RDR graphs, previously identified in \cite{kuzman regular}.
\begin{proposition}\label{basic RDR} The following vertex-transitive graphs are rainbow domination regular:
\begin{enumerate}[$(i)$]
\item For any integer $d\geq 1$, complete bipartite graphs $K_{d,d}$ are $d$-RDR.
\item Cycles $C_n$ are $2$-RDR iff $n\equiv 0\pmod 4$.
\item Prisms $\Prism(n)$ are $3$-RDR iff $n\equiv 0\pmod 6$.
\item M\"obius ladders $\Ml(n)$ are $3$-RDR iff $n\equiv 3\pmod 6$.
\item Wreath graphs $W(n)=C_n[2K_1]$ are $4$-RDR iff $n\geq 4$, $n\equiv 0\pmod 4$.
\end{enumerate}
\end{proposition}
\begin{figure}
\centering\includegraphics[width=15cm]{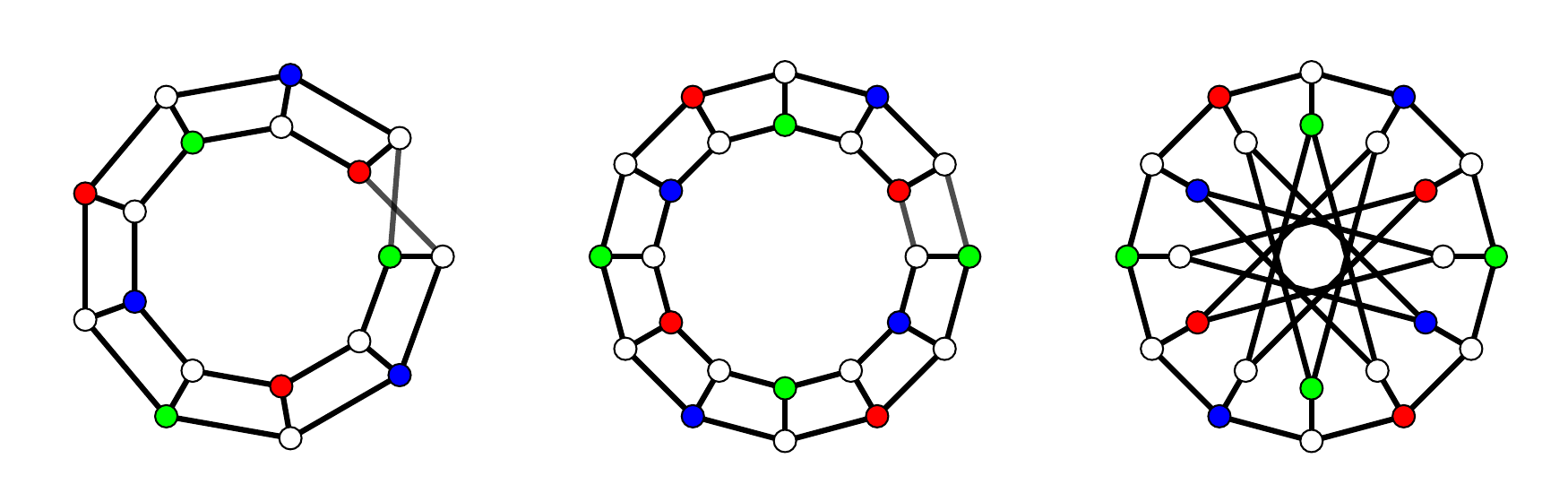}
\caption{M\"obius ladder $\Ml(9)$, prism $\Prism(12)\cong \GP(12,1)$ and Nauru graph $\GP(24,5)$ are examples of vertex-transitive 3-RDR graphs of orders $18$ and $24$.}
\end{figure}

We will now prove two general criteria for a vertex-transitive graph to be $d$-RDR. Let $X$ be a graph and let $H\leq \Aut(X)$ ba a subgroup of automorphisms of $X$. \emph{The quotient graph} $X_H$ of $X$ relative to $H$ is defined as the graph with vertices the orbits of $K$ on $V(X)$ and two orbits adjacent if there is an edge in $E(X)$ with one endpoint in each of these two orbits. 

\begin{theorem}\label{krit1}
Let $X$ be a connected bipartite vertex-transitive graph of order $n$ and degree $d$. If there exists a semiregular subgroup $H\leq Aut(X)$ of order 
${n\over 2d}$ such that the quotient graph $X_H$ is a simple graph, isomorphic to 
$K_{d,d}$, then $X$ is a $d$-RDR graph.
\end{theorem}
\begin{proof}
Suppose this is true. Denote the orbits by $H_1,\ldots,H_d$ and $C_1,\ldots, C_d$, so that every $H_i$ and $C_j$ represent adjacent vertices in $K_{d,d}$. This defines a natural coloring function on $V(X)$ by 
$$f(v)=\begin{cases}\{i\},& v\in C_i;\\
\emptyset,& \text{ otherwise}.
\end{cases}$$ 
Now observe that $H_i\sim C_j$ implies $v\sim v'$ for some $v\in H_i$ and $v'\in C_j$. If $u\in H_i$, then $u=v^h$ for some $h\in H$ and $u$ is adjacent to $v'^h\in C_j$. Thus, any vertex in $H_i$ is adjacent to a vertex in $C_j$, implying that any vertex in $H_i$ has (exactly one) neighbour in $C_j$ for each $j$. By similar argument, the reverse also holds, thus $f$ is a $d$-RDR coloring of $X$.
\end{proof}
\begin{example}
As checked by computation, graph $[36,12]$ from Table~\ref{tab: ord<36} is the smallest vertex-transitive $3$-RDR graph not meeting condition of Theorem~\ref{krit1}, with further examples of orders $48, 60, 72,$ etc. 
On the other hand, graph $[36,11]$ or $X_3$ as defined in subsection~\ref{subs Xn}, is the smallest example of a vertex-transitive $3$-RDR graph, meeting condition of Theorem~\ref{krit1} for two different semiregular subgroups of $\Aut(X)$, thus giving two essentially different $3$-RDR colorings.
\end{example}

\begin{theorem}\label{krit2}
Let $X$ be a connected bipartite vertex transitive graph of order $n$ and degree $d$. Let $V(X)=V_1\cup V_2$  be the bipartition of the vertex set, and  let $G\leq Aut(X)$ be vertex-transitive.
Suppose there exists a block $B\subset V(X)$ for the action of $G$, such that $|B|=\frac{n}{2d}$ and $B\subseteq V_1$, and there is a vertex $v\in B$ such that $N(v)\cap B^g$ is nonempty for exactly $d$ different blocks $B^g$, $g\in G$. Then $X$ is a $d$-RDR graph.
If a block like this exists, it is an orbit of a vertex-transitive subgroup $H\leq G$ of size $|H|\equiv 0\pmod {\frac{n}{2d}}$, containing the vertex-stabilizer $G_v$.
\end{theorem}
\begin{proof}
Observe that for each $g\in G$, the conjugate block $B^g$ is contained in one of the bipartition sets $V_1$, $V_2$. Since there are exactly $2d$ blocks, we denote the blocks for which $N(v)\cap B^g\ne \emptyset$ by $C_1,C_2,\ldots, C_d$, and the remaining blocks by $B=B_1, B_2,\ldots, B_d$.
It is not difficult to see that by above conditions, every vertex $u\in V(X)$ has neighbors in $d$ different blocks. Suppose on the contrary there is a vertex $u$ with neighbors $v_1,v_2$ in the same block $B'$. Since $G$ is vertex-transitive, we have $u^g=v$ for some $g\in V$, and $v_1^g,v_2^g \in B'^g$ are in the same block, contradicting the assumption. The coloring function $f$ on $V$, defined by $f(u)=\{i\}$ if $u \in C_i$ and $\emptyset$ otherwise, is then a $d$-RDR function.
By a well-known result in permutation group theory, any block $B$ for $G$ containing $v$ is an orbit of the setwise stabilizer $H=\{g\in G\colon B^g=B\}$, which also contains $G_v$ and must have size divisible by the order of the block.
\end{proof}

Both criteria are quite efficient to check with existing tools for algebraic computation and available lists of vertex-transitive graphs of small orders and degrees. In fact, VT 3-RDR graphs of orders between 66 and 96 in Table~\ref{tab: bicubic} were identified by applying Theorem~\ref{krit2} to the graphs of \cite{Potocnik Spiga Verret CVT}, as the direct $3$-RDR test by partitioning the vertices was too slow. However, the conditions of Theorem~\ref{krit2} might not be neccessary, so this part of table may be incomplete.

\section{Vertex-transitive 3-RDR graphs}\label{sec VT 3RDR}

As mentioned before, the class of vertex-transitive cubic graphs is well investigated. 
In this section, we identify the vertex-transitive $3$-RDR graphs among the families of the generalized Petersen graphs and honeycomb toroidal graphs, and construct another infinite family of 3-RDR's as Cayley graphs. This will give us partial classification of vertex-transitive $3$-RDR graphs and some open questions.
Before we proceed, we write down a simple observation about the coloring of 6-cycles in a 3-RDR graph.
\begin{lemma}
Let $C$ be a $6$-cycle in a 3-RDR graph. Then any 3-RDR function colors the consecutive vertices of $C$ as 
$\emptyset-\{1\}-\emptyset-\{2\}-\emptyset-\{3\}$ (up to the choice of starting vertex and direction).
\end{lemma} 

\subsection{\sc Generalized Petersen Graphs}
Recall that the \emph{generalized Petersen graphs} $\GP(n,k)$, where $1\leq k\leq n/2$, are cubic graphs of order $2n$ with vertex set $\{u_i,v_i\colon i\in \ZZ_n\}$, and edge set being a disjoint union of outer edges $\{\{u_i,u_{i+1}\colon \in \in\ZZ_n\}$, spokes $\{\{u_i,v_i\}\colon i\in\ZZ_n\}$ and inner edges $\{\{v_i,v_{i+k}\}\colon i\in\ZZ_n\}$.
A simple criteria for a generalized Petersen graph to be 3-RDR was given in \cite{Zerovnik not VT}. We include the proof for completeness.
\begin{theorem}[\v Zerovnik, \cite{Zerovnik not VT}]\label{krit GP} A generalized Petersen graph $\GP(n,k)$, where $1\leq k\leq n/2$, is $3$-RDR iff $n\equiv 0\pmod 6$ and $k\equiv \pm 1\pmod 6$.
\end{theorem}
\begin{proof}
Let $X=\GP(n,k)$ be 3-RDR. Then $X$ is bipartite, thus $n$ is even and $k$ is odd.
Without loss of generality, assume $f(u_i)=f(v_{i+1})=\emptyset$ for $i$ even, and let $f(u_1)=\{1\}$ and $f(u_3)=\{2\}$. Suppose $f(u_5)=\{1\}$.  Then $f(v_2)=f(v_4)=\{3\}$. This implies that $f(u_{2\pm k}),f(u_{4\pm k})\ne \{3\}$, forcing that $f(v_{3\pm k})=\{3\}$, a contradiction since both are neighbors of $v_3$. Thus $f(u_5)=\{3\}$ and by symmetry, the outer cycle must be colored by pattern $\emptyset-\{1\}-\emptyset-\{2\}-\emptyset-\{3\}$, hence $n\equiv 0\pmod 6$. Also, coloring of the outer cycle now implies that $f(v_0)=\{2\}$, $f(v_2)=\{3\}$, $f(v_4)=\{1\}$, etc, so the only possible coloring function is
$$f(u_i)=f(v_{i+3})=\begin{cases}
\emptyset, &i\equiv 0\pmod 2,\\
\{1\},&i\equiv 1\pmod 6,\\
\{2\},&i\equiv 3\pmod 6,\\
\{3\},&i\equiv 5\pmod 6.
\end{cases}$$
If $k\equiv 3\pmod 6$, then non-colored vertex $v_1$ has neighbors $u_1$ and $v_4$ of the same color. Hence we must have $k\equiv\pm 1\pmod 6$. Conversely, for $n\equiv 0\pmod 6$ and $k\equiv\pm 1\pmod 6$, the above coloring function is 3-rainbow dominating of weight $|V(X)|/2$. This shows that $X$ is 3-RDR.
\end{proof}

As noted in \cite{Zerovnik not VT}, it follows from this result that non-vertex-transitive $3$-RDR graphs exist, with $\GP(18,5)$ and $\GP(18,7)$ of order $36$ being the smallest of this kind. In fact, one could explicitly describe some infinite families of non-vertex-transitive 3-RDR graphs by parameters of $\GP(n,k)$. However, we pursue to classify vertex-transitive 3-RDR graphs among the generalized Petersen graphs.

Recall that the \emph{girth} of graph, denoted by $\girth(X)$, is the length of its shortest cycle. Such cycles are called \emph{girth cycles}. For any edge $e$, we denote by $\eps(e)$ the number of girth cycles containing $e$. Following \cite{Potocnik Vidali girth reg} and \cite{Potocnik Vidali girth six}, the \emph{signature of the vertex} $v$ is the sequence $\signat(v)=(\eps(e_1),\eps(e_2),\ldots,\eps(e_k))$, where $e_i$ are edges, incident with $v$ and ordered in such a manner that $\eps(e_1)\leq \eps(e_2)\leq \ldots\leq \eps(e_k)$. If graph is regular and all vertex signatures are equal, the graph is called \emph{girth regular}. In such case, the vertex signature is called the \emph{signature of the graph} $X$, which we will denote by $\signat(X)$. Obviously, every vertex-transitive graph is girth regular, and its signature is a natural invariant, describing certain combinatorial properties of the graph.

We now combine this with some well-known results on generalized Petersen graphs to obtain a classification of vertex-transitive 3-RDR generalized Petersen graphs by their girth and signature. Recall that generalized Petersen graphs are vertex-transitive iff $k^2\equiv \pm 1\pmod n$ or $(n,k)=(10,2)$ (see \cite{Nedela GP}), and note that the girths of generalized Petersen graphs were determined in \cite{Ferrero Hanusch 2014}. 

\begin{theorem}\label{3RDR GP}\label{thm: 3RDR VT GP}
Let $X\cong \GP(n,k)$ be a generalized Petersen graph, such that $X$ is 3-RDR and vertex-transitive.
Then $X$ is a Cayley graph and one of the following holds:
\begin{enumerate}[(i)]
\item $\girth(X)=4$, $\signat(X)=(1,1,2)$ and $X\cong \GP(n,1)\cong \Prism(n)$ with $n\equiv 0\pmod 6$.
\item $\girth(X)=6$, $\signat(X)=(2,2,2)$ and $X\cong \GP(n, {n\over 2}-1)$ with $n\equiv 0\pmod {12}$.
\item $\girth(X)=8$, $\signat(X)=(8,8,8)$ and $X\cong \GP(24,5)$.
\item $\girth(X)=8$, $\signat(X)=(10,11,11)$ and $X\cong \GP(24,7)$.
\item $\girth(X)=8$, $\signat(X)=(2,2,4)$ and $X\cong \GP(n,k)$ with $n\equiv 0\pmod {18}$, $n\geq 72$, $17\leq k\leq {n\over 2}-2$, $k\equiv\pm 1\pmod {18}$ and $k^2\equiv 1\pmod n$.
\item $\girth(X)=8$, $\signat(X)=(5,5,6)$ and $X\cong \GP(n,k)$ with $n\equiv \pm 6 \pmod {18}$, $n\geq 30$, $5\leq k\leq {n\over 2}-2$, $k\equiv\pm 1\pmod 6$ and $k^2\equiv 1\pmod n$.

\end{enumerate}

\end{theorem}
\begin{proof}
Let $X=\GP(n,k)$ be 3-RDR and vertex transitive. Then $n=6m$ and $k\equiv \pm 1\pmod 6$.
Recall first that by \cite{Nedela GP}, graph $\GP(n,k)$ is vertex-transitive iff $k^2\equiv \pm 1\pmod n$ or $(n,k)=(10,2)$, and is a Cayley graphs if $k^2\equiv 1\pmod n$. Suppose that $X$ is not a Cayley graph. Then $k^2\equiv -1\pmod n$. Since $k=6\ell\pm 1$, this implies $k^2=36\ell^2\pm12\ell+2\equiv 0\pmod {6m}$, a contradiction. Therefore, we can assume that $k^2\equiv 1\pmod{ 6m}$ and $X$ is a Cayley graph.
By results of \cite{Ferrero Hanusch 2014}, generalized Petersen graphs $\GP(n,k)$ have girth $\leq 8$.
Since in our case $X$ is bipartite, the girth of $X$ is necessarily even. 
\begin{itemize}
\item If $\girth(X)=4$, by \cite{Ferrero Hanusch 2014}, we either have $k=1$, or $n=4k$. The latter case implies a contradiction with $k^2\equiv 1\pmod n$, as this implies $k^2- 4\ell k-1=0$ for some integers $k,\ell$, but solving for $k$ we obtain non-integral solutions. Hence $k=1$ is the only possibility and the graph obtained is a prism. It is easy to see that the signature is $(1,1,2)$ for $n\ne 4$.
\item If $\girth(X)=6$, by \cite{Ferrero Hanusch 2014}, we either have $k=3$ or $n=6k$ or $n=2k+2$. If $k=3$, then $k^2=9\equiv 1\pmod {6n}$, a contradiction. If $n=6k$, then $k^2\equiv 1\pmod n$ has no solutions. The last case $n=2k+2$ implies $n=2k+2=6m$, hence $k=3m-1$, so $m$ must be even. In this case, we have that $k^2=9m^2-6m+1\equiv 1\pmod {6m}$, complying with vertex-transitivity of $X$. In this case, we can count $6$-cycles containing a given edge to see that $(2,2,2)$ is the signature of the graph.

Generally, there are 3 possible types of $6$-cycles an a generalized Petersen graph, containing 1, 2 or 3 consecutive outer edges. Let $C$ be a $6$-cycle in $\GP(6m, 3m-1)$.
\begin{figure}[ht!]
\centering
\includegraphics[width=5cm]{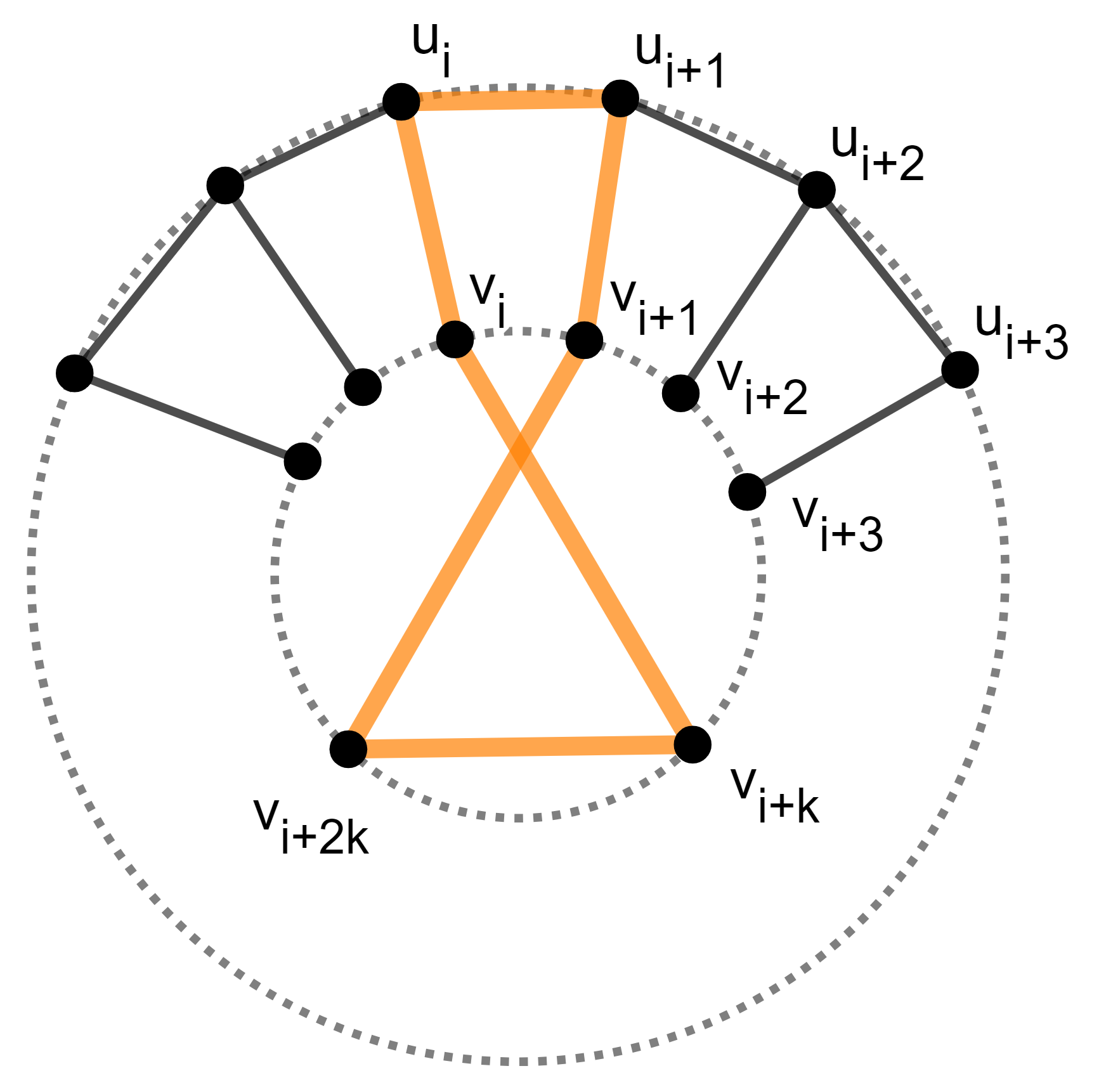}
\includegraphics[width=5cm]{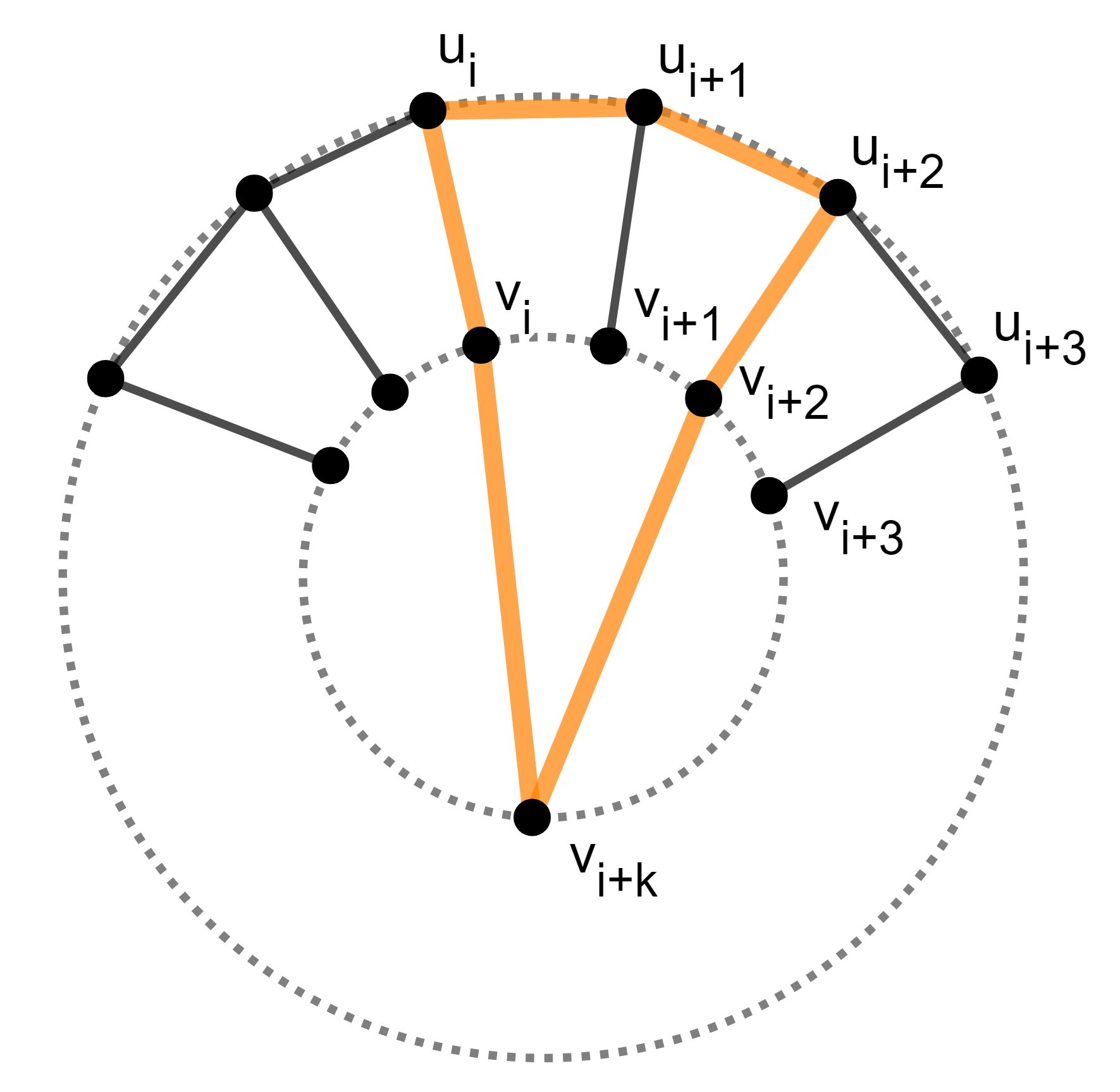}
\includegraphics[width=5cm]{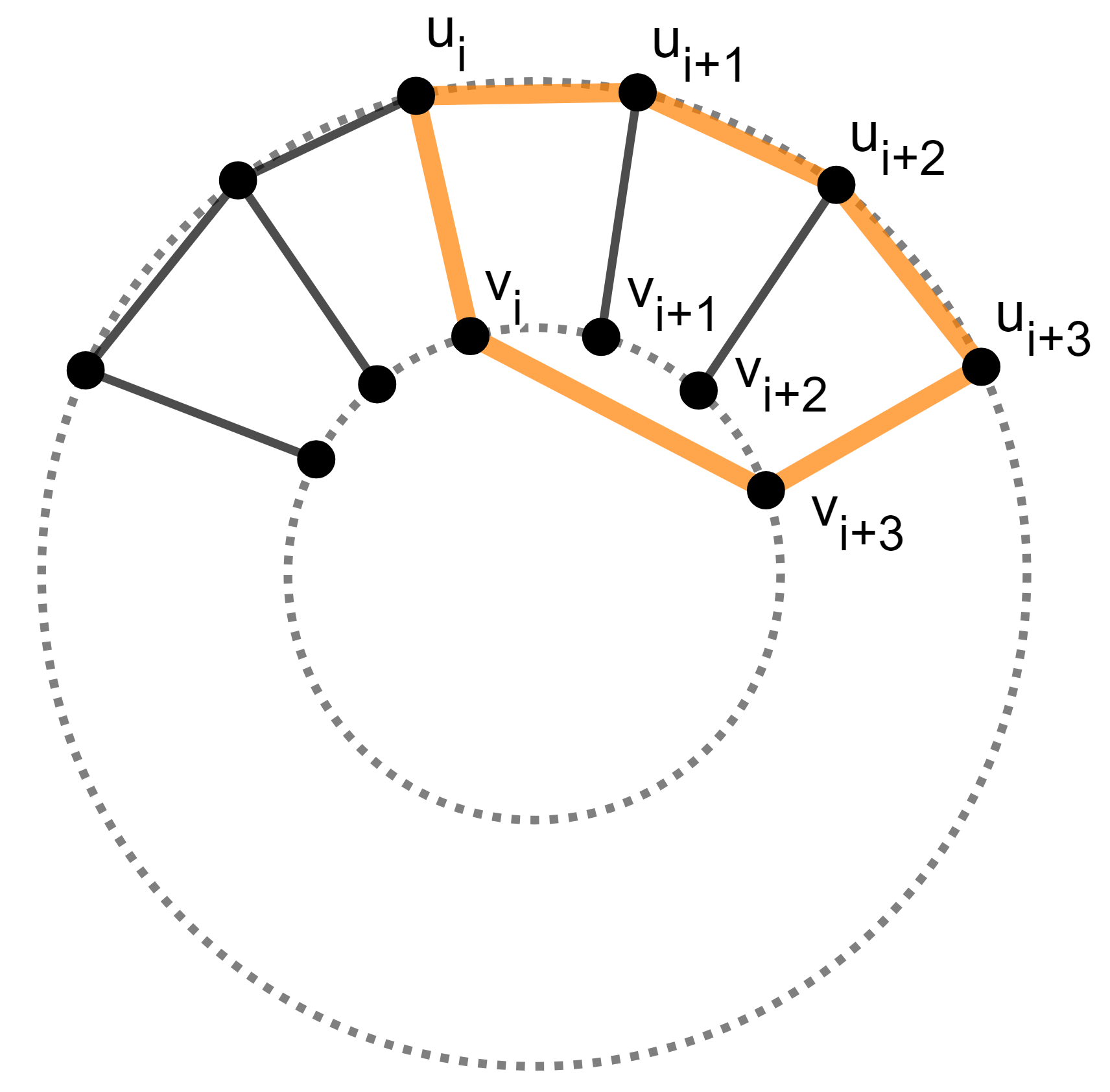}
\caption{The three types of $6$-cycles in generalized Petersen graphs}\label{GP 6c}
\end{figure}

If $C$ contains $1$ outer edge, then it contains 2 spokes and 3 consecutive inner edges, so the sequence of vertices is $u_iv_iv_{i+k}v_{i+2k}v_{i\pm 1}u_{i\pm1}$ (Figure~\ref{GP 6c}, left). Equality $v_{i+3k}=v_{i\pm 1}$ implies $3k\equiv \pm 1\pmod {6m}$ and since $k^2\equiv 1\pmod {6m}$, by squaring we obtain $8\equiv 0\pmod {6m}$, a contradiction.

If $C$ contains $2$ outer edges, they must be consecutive and $C$ also contains 2 spokes and $2$ inner edges (Figure~\ref{GP 6c} in the middle). The sequence of vertices $u_iv_i v_{i+k} v_{i\pm 2} u_{i\pm 2}u_{i\pm 1}$ and equality $v_{i+2k}=v_{i\pm 2}$ implies $2k\equiv \pm 2\pmod {6m}$, hence $k\equiv \pm 1\pmod {3m}$ and thus $k=3m-1$, which is fulfilled in our case. By symmetry, every edge is contained in exactly two $6$-cycles of this type.

If $C$ contains $3$ outer edges, by similar reasoning the implied equation $k\equiv \pm 3\pmod {6m}$ leads to a contradiction by squaring.

Hence, each edge in $\GP(6m,3m-1)$ is contained in exactly $2$ different $6$-cycles, and three edges meeting in a vertex determine the signature of graph to be $(2,2,2)$.

\item If $\girth(X)=8$, we study possible $8$-cycles in $X$ in similar fashion. Let $C$ be an $8$-cycle in $\GP(6m,k)$. First we show that with our restrictions on $n,k$, some types of $8$-cycles are possible only for certain values of $n$ or not possible at all.

\begin{figure}[h!]
\centering
\includegraphics[width=5cm]{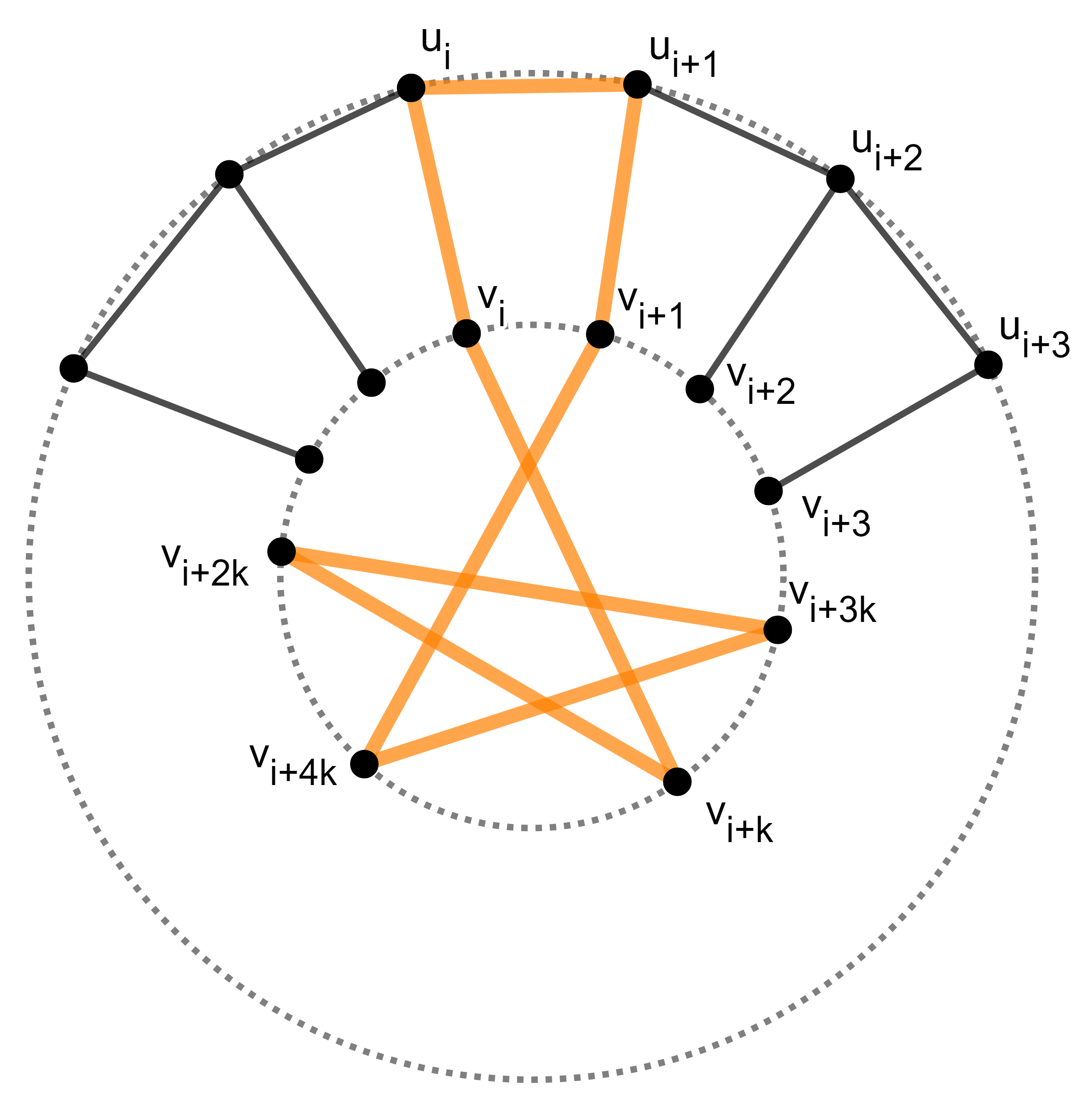}
\includegraphics[width=5cm]{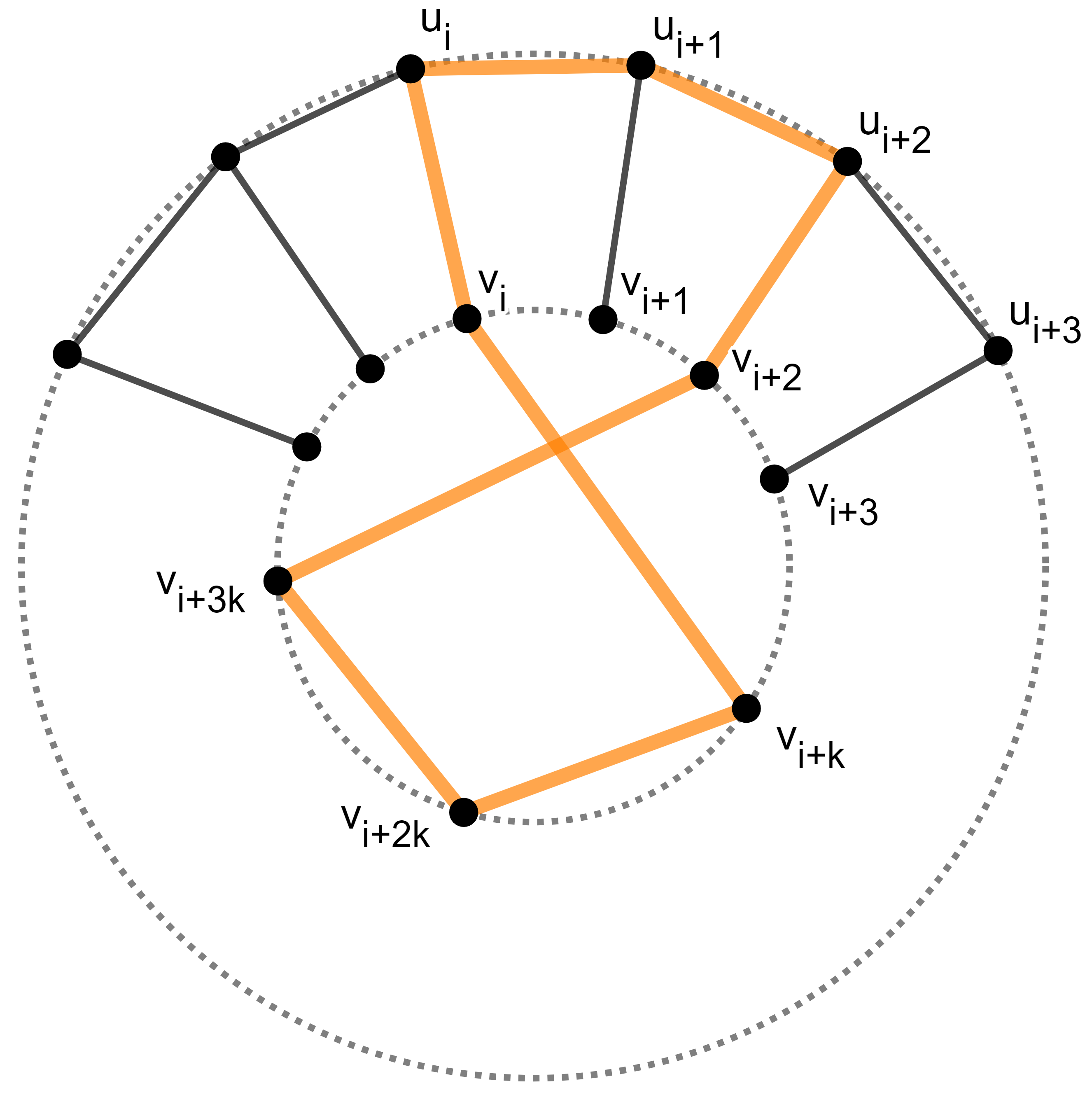}
\includegraphics[width=5cm]{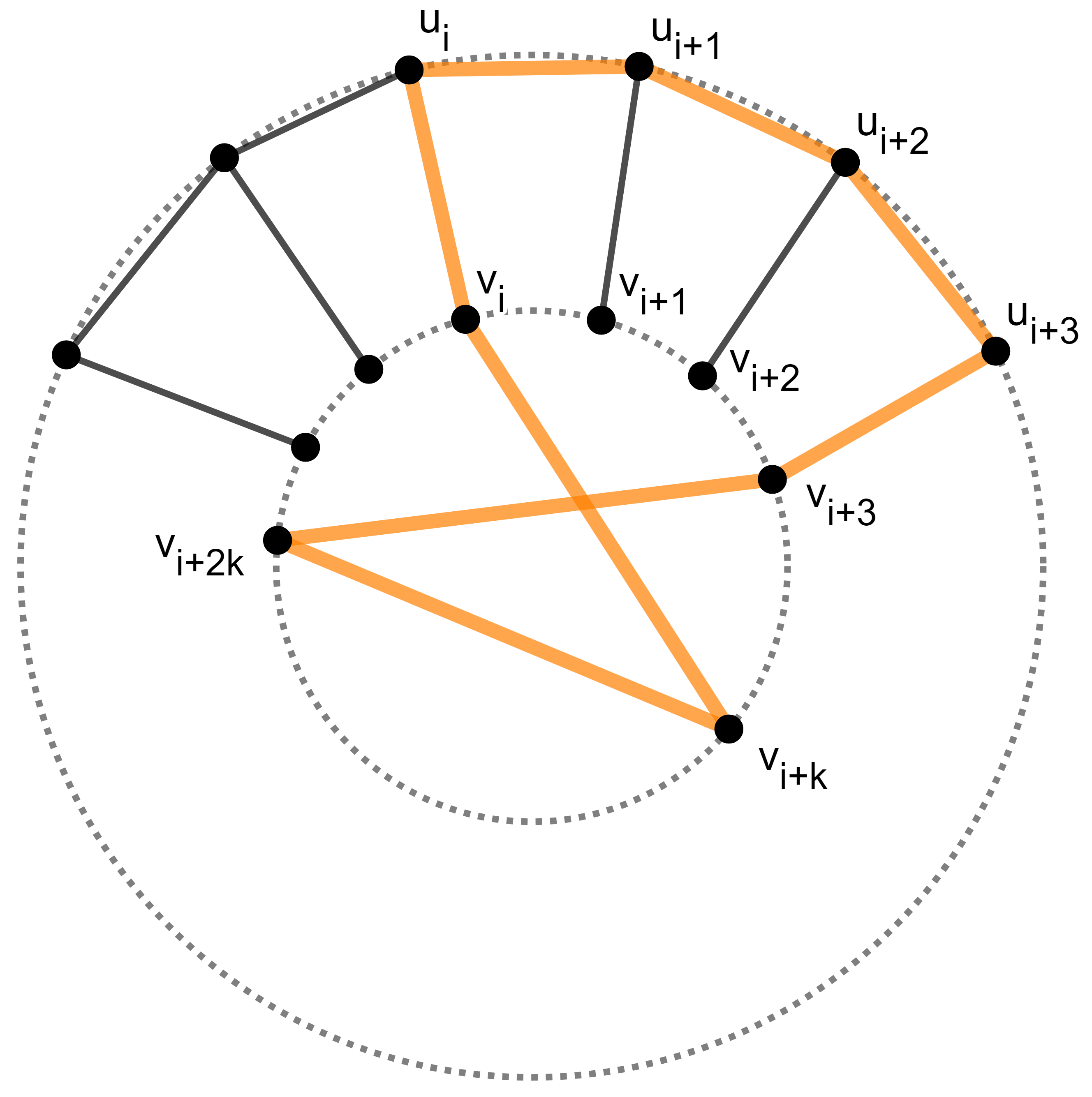}
\caption{$8$-cycles with exactly $1$, or $2$ or $3$ consecutive outer edges in generalized Petersen graphs}\label{GP 8ca}
\end{figure}
If $C$ contains 0 outer edges, then it contains $8$ inner edges and no spokes. This implies $8k\equiv 0\pmod{6m}$ and hence $64k^2\equiv 0\pmod{6m}$ while $k^2\equiv 1\pmod{6m}$. It follows that $6m|64$ and hence $3|32$, a contradiction.

If $C$ contains exactly 1 outer edge, then it also contains 2 spokes and 5 consecutive inner edges (see Figure~\ref{GP 8ca} on the left), implying that $5k\equiv \pm 1\pmod{6m}$. It follows that $25k^2\equiv 1\pmod {6m}$, so $25\equiv 1\pmod {6m}$ and hence $m|4$, implying that $m=1,2$ or $4$. However, since $k\ne 1, 3m-1$, the only possible graphs of girth $8$ with cycles of this type are $\GP(24,5)$ and $\GP(24,7)$. 

If $C$ contains 2 consecutive outer edges, then it contains 2 spokes and 4 consecutive inner edges (Figure~\ref{GP 8ca} in the middle). In similar fashion as before, we get $4k\equiv\pm 2\pmod{6m}$ and hence $16k^2\equiv 16\equiv 4\pmod{6m}$. Then $m|2$, again contradicting that $\girth(X)=8$.
 
If $C$ contains $3$ consecutive outer edges, then it contains $2$ spokes and $3$ consecutive inner edges. It follows that $3k\equiv \pm 3\pmod{6m}$. In this case, using $k^2\equiv 1\pmod{6m}$ yields no contradiction, but we can reduce the earlier equation to $k\equiv \pm 1\pmod {2m}$. Since $1<k< 3m-1$, two possibilities for $k$ are obtained: $k=2m\pm 1$. Now 
$k^2=4m^2\pm 4m+1\equiv 1 \pmod{6m}$
implies that $6m|4m^2+4m$ and so $3|(m\pm 1)$. Cycles of this type are therefore possible for $m\equiv \pm 1\pmod 3$ (and hence $n\equiv \pm 6\pmod{18}$), but not for $m\equiv 0\pmod 3$. We also point out that in in the latter case, we have $n\equiv 0\pmod 18$, and from $n|k^2-1$ we obtain $k\equiv \pm 1\pmod{18}$ by a short computation. 

By similar computations, cases when $C$ contains $4$ or $5$ consecutive outer edges are only possible for $m|2$ or $m|4$, the latter case applies to graphs $\GP(24,5)$ and $\GP(24,7)$ of girth $8$. Moreover, $8$-cycles with $6$ or more consecutive outer edges are also impossible for obvious combinatorial reasons.

However, it is possible that $C$ contains two non-consecutive outer edges. In such case, the structure of $8$-cycle $C$ is necessarily \emph{outer-spoke-inner-spoke-outer-spoke-inner-spoke}, and we have two further possibilities (see Figure~\ref{GP 8cb}). 

The sequence of vertices $u_iv_iv_{i+k}u_{i+k}u_{i+k-1}v_{i+k-1}v_{i-1}u_{i-1}$ is an $8$-cycle in $\GP(6m,k)$ (for any $k$, $1\leq k \leq 3m$, see Figure~\ref{GP 8cb} on the left). Note that, by symmetry, each spoke of the graph is contained in 4 different $8$-cycles of this type, and each outer or inner edge is contained in 2 different $8$-cycles of this type.

Finally, another possible type of $8$-cycle is given by the sequence
$u_iv_iv_{i+k}u_{i+k}u_{i+k+1}v_{i+k+1}v_{i+2k+1}u_{i+2k+1}$ (Figure~\ref{GP 8cb} on the right), but then $k=3m-1$ and the girth of graph is $6$. 
\begin{figure}[h!]
\centering
\includegraphics[width=5cm]{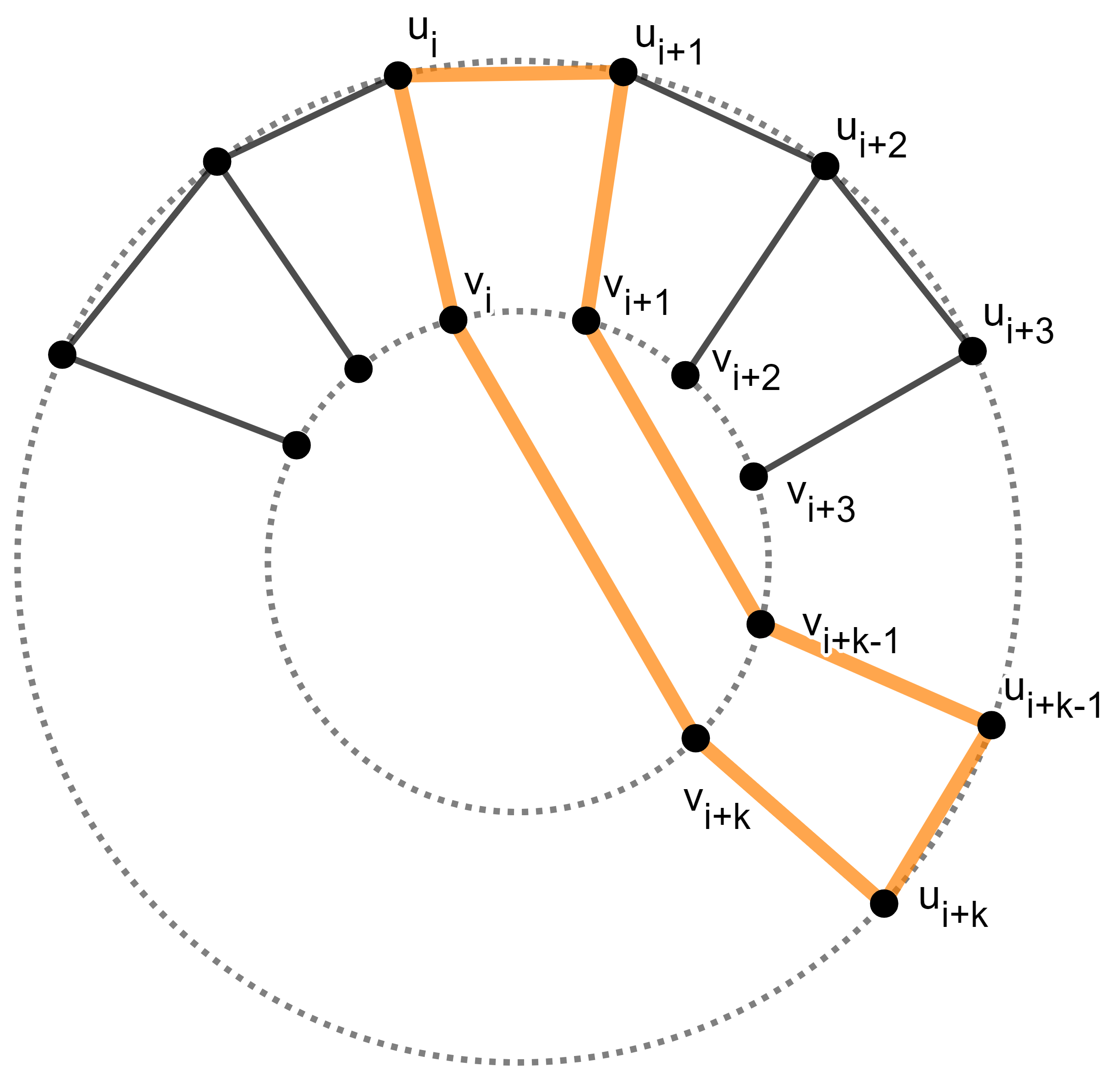}
\includegraphics[width=5cm]{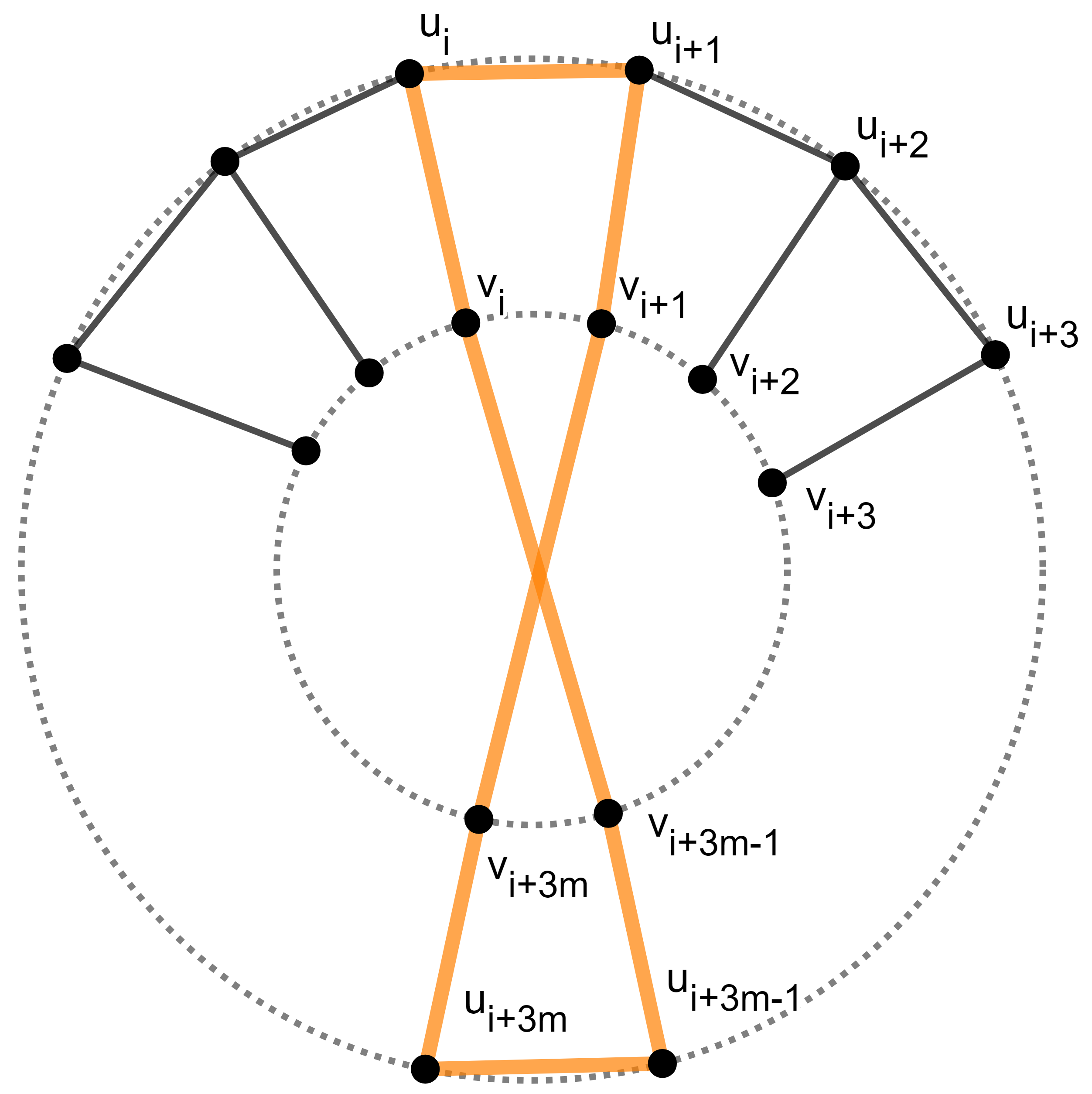}
\caption{$8$-cycles with exactly $2$ nonconsecutive outer edges in generalized Petersen graphs}\label{GP 8cb}
\end{figure}
Combining and counting above solutions for girth 8, we have see that for $m=4$, the signature is $(8,8,8)$ for $\GP(24,5)$ and $(10,11,11)$ for $\GP(24,7)$,  while for $m\geq 5$, we have $(5,5,6)$ if $m\equiv \pm 1\pmod 3$ and $(2,2,4)$ for $m\equiv 0\pmod 3$.
Moreover, the smallest appropriate pair $(n,k)$ is $(30,11)$ in the first case, and $(72,17)$ in the second case. 
 \end{itemize}
\end{proof}

\begin{corollary}
For girths $4,6,8$, there are infinitely many vertex-transitive 3-RDR graphs of each girth.
\end{corollary}
Next, we proceed to construct and classify more graphs of girth 6.

\subsection{\sc Honeycomb toroidal graphs}
The honeycomb toroidal graphs or HTG's are a specific family of vertex-transitive cubic graphs, defined in \cite{Sparl} as follows. 
Let $m,n,\ell$ be integers, such that $m\geq 1$, $n\geq 4$ is even, and $0\leq \ell\leq n/2$, where $\ell\equiv m\pmod 2$ (so $\ell$ and $m$ have same parity). The \emph{honeycomb toroidal graph }$\HTG(m,n,\ell)$ has vertex set $V=\{v_{i,j}\mid i\in {\mathbb Z}_m, j \in {\mathbb Z}_n\}$ and edge set $E$ being the disjoint union of 
\begin{itemize}
\item \emph{vertical edges}
$\{\{ v_{i,j},v_{i,j+1}\}\mid  i\in {\mathbb Z}_m, j \in {\mathbb Z}_n\}$,
\item \emph{horizontal edges} $\{\{v_{i,j},v_{i+1,j}\} \mid  i\in{\mathbb Z}_m\setminus\{m-1\}, j\in{\mathbb Z}_n, i+1\equiv j \pmod 2\}$, and
\item  \emph{jumps} 
$\{\{v_{m-1,j},v_{0,j+\ell}\} \mid j \in{\mathbb Z}_n, j\equiv m\pmod 2\}$.
\end{itemize}
It is known that all HTG graphs are bipartite and cubic, vertex-transitive and can be described as Cayley graphs of generalized dihedral groups, with respect to a connection set of three involutions. Moreover, they have girth $\leq 6$ and can be embedded on the torus with hexagons as their faces (see \cite{Alspach 1}, \cite{Alspach 2}). 

\begin{theorem}\label{3RDR HTG}\label{krit HTG}
Let $X=\HTG(m,n,\ell)$. Then $X$ is a 3-RDR graph iff 
$m$ is even and $n\equiv \ell \equiv 0\pmod 6$, or $m$ is odd, $n\equiv0\pmod 6$ and $\ell\equiv 3 \pmod 6$.
\end{theorem}
\begin{proof}
First, let $m=1$. Suppose that $X$ is a 3-RDR graph. Then $n=|V(X)|\equiv 0\pmod 6$ and $n\geq 6$.
Note that $X$ contains only vertical edges $\{v_{0,j},v_{0,j+1}\}$ and jumps $\{v_{0,j},v_{0,j+\ell}\}$. Wlog suppose that $f(v_{0,j})=\emptyset$ for $j$ even, and $f(v_{0,1})=\{1\}$, $f(v_{0,3})=\{2\}$. Then $f(v_{0,2+\ell})=\{3\}$, since $v_{0,2}$ must have neighbours of all colors. If $f(v_{0,5})=\{1\}$, then also $f(v_{0,4+\ell})=\{3\}$, a contradiction as then $v_{0,3+\ell}$ has two neighbors of same color. This implies $f(v_{0,5})=\{3\}$. Repeating the argument, we see that the vertical cycle is colored by a pattern 
$\emptyset-\{1\}-\emptyset-\{2\}-\emptyset-\{3\}-\cdots$, that is, 
$$f(v_{0,j})=\begin{cases} 
\>\>\emptyset,&j\equiv 0\pmod 2,\\
\{1\},& j\equiv 1\pmod 6,\\
\{2\},& j\equiv 3\pmod 6,\\
\{3\},& j\equiv 5\pmod 6.
\end{cases}$$
By definition, $\ell$ is odd. If $\ell\equiv 1\pmod 6$, then vertex $v_{0,0}$ has two neighbours of same color as $f(v_{0,1})=f(v_{0,\ell})=\{1\}$, a contradiction. Case $\ell\equiv 5\pmod 6$ is similar, implying that condition $\ell\equiv 3\pmod 6$ is necessary. For the reverse implication, it is easy to see that conditions $n\equiv 0\pmod 6$ and $\ell\equiv 3\pmod 6$ imply the above coloring function is 3RD of weight $|V(X)|/2$, hence $X$ is 3-RDR.

\begin{figure}[ht!]
\centering
\includegraphics[height=5cm]{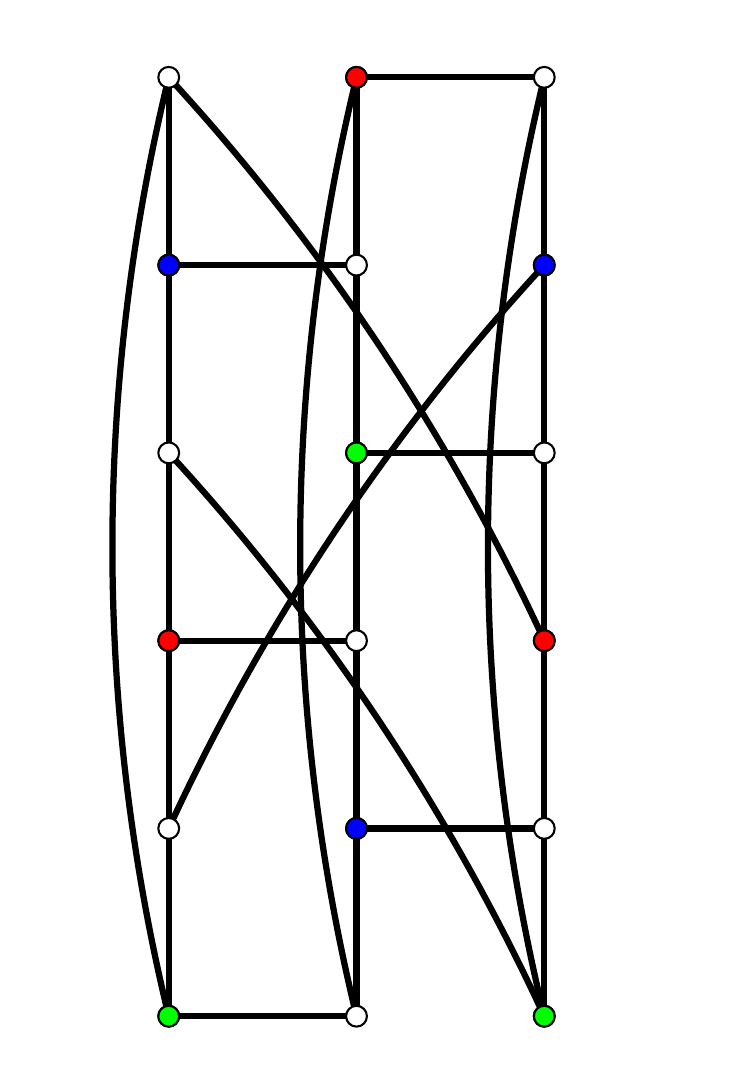}~
\includegraphics[height=5cm]{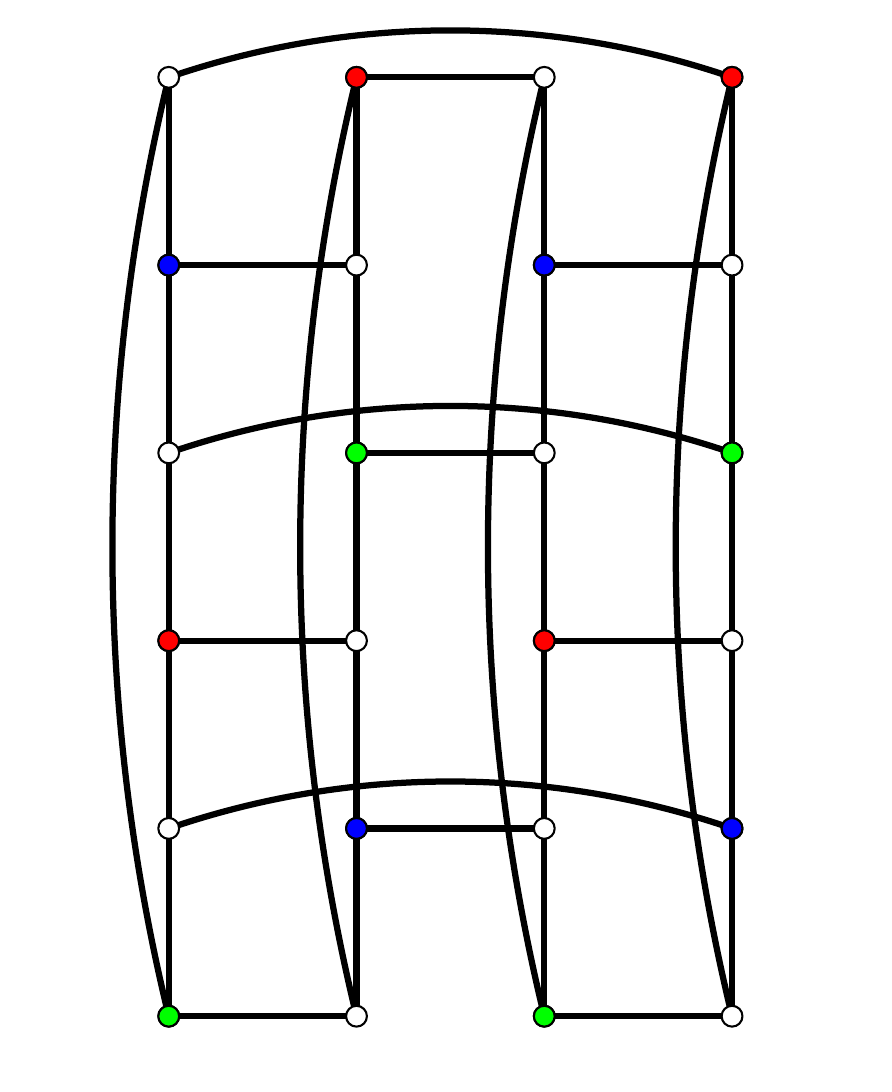}~
\includegraphics[height=5cm]{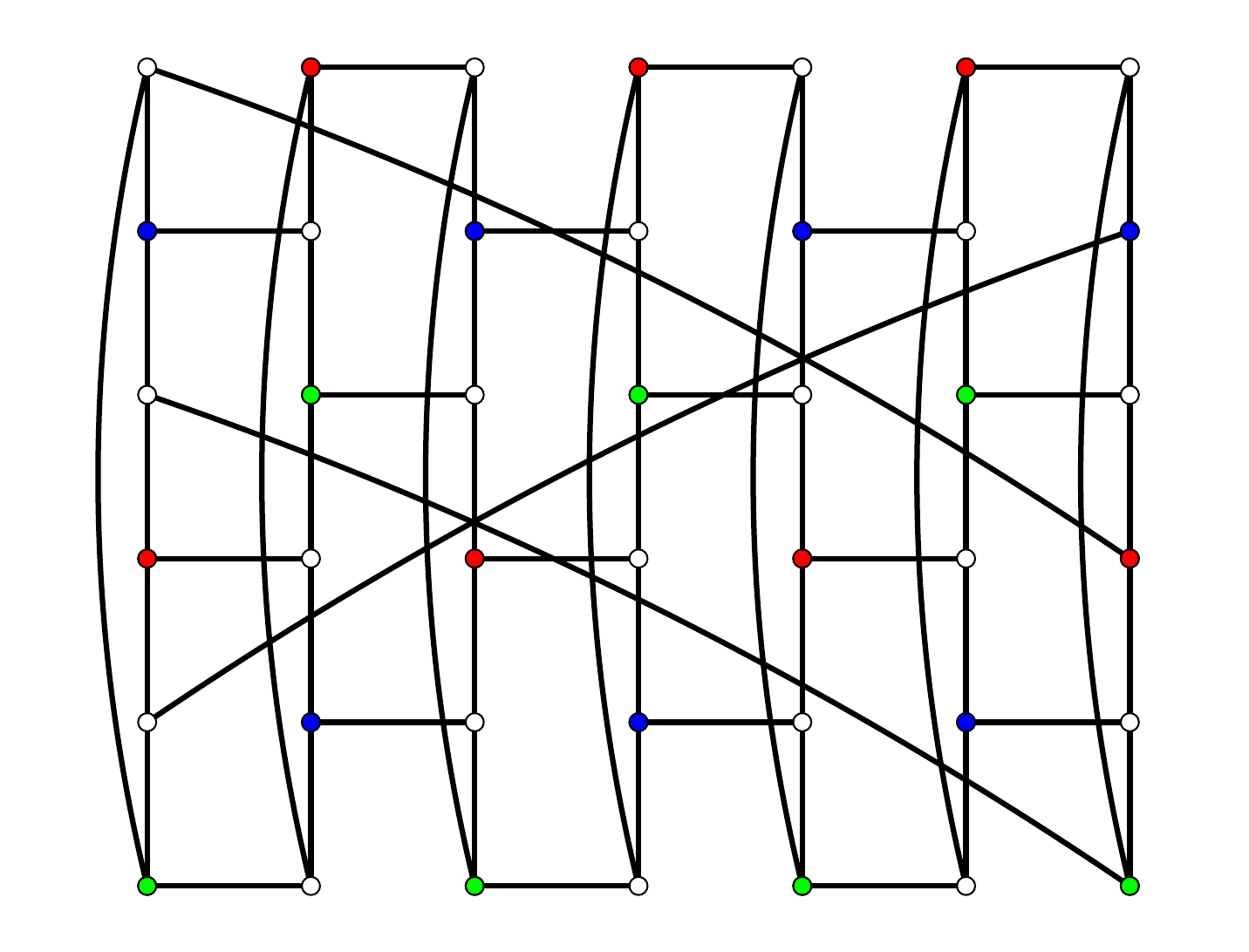}\caption{3-RDR colorings of graphs $\HTG(3,6,3)$, $\HTG(4,6,0)$ and $\HTG(7,6,3)$.}
\label{fig: 8cycles1}\label{fig HTG}
\end{figure}

Now let $m\geq 2$. Wlog suppose that $f(v_{i,j})=\emptyset$ for $j+i$ even and let $f(v_{0,1})=\{1\}$, $f(v_{0,3})=\{2\}$. Then $f(v_{1,2})=\{3\}$, since $v_{0,2}$ must have neighbours of all colors. Repeating the argument, we obtain $f(v_{1,0})=\{2\}$, $f(v_{1,4})=\{1\}$ and $f(v_{0,5})=\{3\}$ (note that we must have $n>4$). Again, we see that the vertical cycles must be colored by a pattern 
$\emptyset-\{1\}-\emptyset-\{2\}-\emptyset-\{3\}-\cdots$, and hence $n\equiv 0\pmod 6$. However, the colors in two consecutive columns $i$ and $i+1$ ($i=0,\ldots, n-2$) are shifted by three steps, so the appropriate coloring function is
$$f(v_{i,j})=\begin{cases} 
\>\>\emptyset,&i+j\equiv 0\pmod 2,\\
\{1\},& i \text{ even and } j\equiv 1\pmod 6, \text{ or $i$ odd and } j\equiv 4\pmod 6,\\
\{2\},&  i \text{ even and } j\equiv 3\pmod 6, \text{ or $i$ odd and } j\equiv 0\pmod 6,\\
\{3\},& i \text{ even and } j\equiv 5\pmod 6, \text{ or $i$ odd and } j\equiv 2\pmod 6.
\end{cases}$$
In order to obtain conditions on $\ell$, we now inspect jumps for cases $m$ even and odd separately. 
For $m$ even, $\ell$ must also be even. Suppose $\ell\equiv 2\pmod 6$. Then $f(v_{m-1,0})=\{2\}=f(v_{0,1+\ell})$, so two neighbors of $v_{0,\ell}$ have the same color, a contradiction. With similar arguments, we see that $\ell\equiv 0\pmod 6$ for $m$ even, and $\ell\equiv 3\pmod 6$ for $m$ odd. 
\end{proof}
We now sharpen this result by further classifying the 3-RDR HTG graphs by their girth and signature. Also, we can identify the generalized Petersen graphs among them.
\begin{theorem}\label{3RDR HTG sgn}\label{thm HTG}
Let $X=\HTG(m,n,\ell)$ be a 3-RDR graph. Then $m$ is even and $n\equiv \ell \equiv 0\pmod 6$, or $m$ is odd, $n\equiv0\pmod 6$ and $\ell\equiv 3 \pmod 6$. Moreover, one of the following holds:
\begin{enumerate}[(i)]

\item $\girth(X)=4$, 
$\signat(X)=(1,1,2)$, and $X$ is 
isomorphic to a prism or a M\"obius ladder,
$$X\cong \HTG(1,2n,3)\cong
\begin{cases}
\Prism(n)\cong \GP(n,1)\cong \HTG(2,n,0),& n\equiv 0\pmod 6,\\
\Ml(n)\cong \HTG(1,2n,n),&n\equiv 3\pmod 6.
\end{cases}
$$
\item $\girth(X)=6$, 
$\signat(X)=(4,4,4)$ and $X\cong \HTG(3,6,3)$, also known as the Pappus graph.

\item $\girth(X)=6$, $\signat(X)=(3,3,3)$ and 
$$X\cong \HTG(3,n,3)\cong
\begin{cases}
\HTG(\frac{n}{2},6,0),& n\equiv 0\pmod {12},\\
\HTG(\frac{n}{2},6,3),& n\equiv 6\pmod {12},\\
\end{cases}
$$ where $n>6$.
\item $\girth(X)=6$, 
$\signat(X)=(2,3,3)$, and $X$ is isomorphic to
$$X\cong \HTG(m,6,0)\cong 
\begin{cases}
\HTG(1,6m,2m-1),&m\equiv 2\pmod 6,\\
\HTG(1,6m,2m+1),&m\equiv 4\pmod 6,\\
\end{cases}$$
where $m$ is even, or 
$$X\cong \HTG(m,6,3)\cong 
\begin{cases}
\HTG(1,6m,2m-1),&m\equiv 5\pmod 6,\\
\HTG(1,6m,2m+1),&m\equiv 1\pmod 6,\\
\end{cases}$$
where $m\geq 5$ is odd.
\item $\girth(X)=6$, $\signat(X)=(2,2,2)$ and $X$ is isomorphic to 
the generalized Petersen graph 
$$X\cong\HTG(2,n,\frac{n}{2})\cong\GP(n,\frac{n}{2}-1),$$
where $n\equiv 0 \pmod 4,$
or $X\cong\HTG(m,n,\ell)$ has parameters $m,n,\ell$ different from any above and is not isomorphic to any generalized Petersen graph.
\end{enumerate}
\end{theorem}
\begin{proof}
Observe that the basic properties of parameters $m,n,\ell$ come from Theorem \ref{3RDR HTG}, and suppose first that $\girth(X)=4$. 
By \cite{Alspach 1}, HTG graphs of girth $4$ are classified as 
prisms, $\Prism(n)\cong \HTG(2,n,0)\cong \HTG(1,2n,3)$, where $n\geq 4$ is even, M\"obius ladders  
 $\Ml(n)\cong \HTG(1,2n,n)\cong \HTG(1,2n,3)$, where $n\geq 3$ is odd,
 or generalized prisms $\GPr(m)\cong \HTG(m,4,\ell)\cong \HTG(2,2m,2)\cong \HTG(1,4m,2m-1)$, where $m\geq 2$, and $\ell\in \{0,2\}$ for $m$ even and $\ell=1$ for $m$ odd. We can read from the parameters or check directly that the generalized prisms are never 3-RDR. The result for prisms and M\"obius ladders follows from \cite{kuzman regular}, and it is easy to see that $\signat(X)=(1,1,2)$ in this case.
\medskip

Now suppose that $\girth(X)=6$. In order to determine the possible signatures, we analyze the possible $6$-cycles in a 3-RDR HTG graph (refer also to Figure~\ref{fig HTG} for illustration of different cases). 
\begin{itemize}
\item A 6-cycle $C$ with $6$ vertical edges is only possible if $n=6$. \item A 6-cycle $C$ with $5$ vertical edges is not possible.
\item A 6-cycle $C$ with $4$ vertical edges can contain $2$ horizontal edges or $2$ jumps in the sequence up-up-right-down-down-left, hence  $C=(v_{i,j}v_{i,j+1}v_{i,j+2}v_{i+1,j+2}v_{i+1,j+1}v_{i+1,j})$ or\\
$C=(v_{m-1,j}v_{m-1,j+1}v_{m-1,j+2}v_{0,j+2+\ell}v_{0,j+1+\ell}v_{0,j+\ell})$.
\item A 6-cycle $C$ with $3$ vertical edges is also possible with 2 horizontal edges and a jump in the sequence up-right-up-right-up-jump. Observe that in this case, we must necessarily have $m=3$ and $\ell=3$.
\item Other $6$-cycles are not possible, as each horizontal edge or jump must be followed by a vertical edge.
\end{itemize}
It is now easy to see that each vertical edge is contained in at most $4$ different $6$-cycles, and this only happens when $n=6$, $m=\ell=3$.
In this case, $X=\HTG(3,6,3)$ is the Pappus graph with signature $\signat(X)=(4,4,4)$ (note that $X$ is also arc-transitive).

If $n>6$ and $m=\ell=3$, then each vertical edge is contained in $3$ different $6$-cycles, and the same for horizontal edges and jumps. This gives us graphs $X\cong\HTG(3,n,3)$ with $n\equiv 0\pmod 6$ and $\signat(X)=(3,3,3)$. In this case, we can also show that $X$ is isomorphic to $\HTG(n/2,6,0)$ if $n/2$ is even, or to $\HTG(n/2,6,3)$ if $n/2$ is odd. Let $n/2$ be even, then $n\equiv 0\pmod{12}$.
In this case, we can define a graph isomorphism $\varphi\colon \HTG(3,n,3)\to \HTG(n/2,6,0)$ by
$$\varphi(v_{0,j})=
\begin{cases}
v_{\frac{j}{2},\frac{j}{2}}, &j\equiv 0\pmod 2,\\
v_{\frac{j-1}{2},\frac{j+1}{2}}, &j\equiv 1\pmod 2,\\
\end{cases}
\varphi(v_{1,j})=
\begin{cases}
v_{\frac{j-2}{2},\frac{j+4}{2}}, &j\equiv 0\pmod 2,\\
v_{\frac{j-1}{2},\frac{j+3}{2}}, &j\equiv 1\pmod 2,\\
\end{cases}
$$
and
$$
\varphi(v_{2,j})=
\begin{cases}
v_{\frac{j-2}{2},\frac{j+6}{2}}, &j\equiv 0\pmod 2,\\
v_{\frac{j-3}{2},\frac{j+7}{2}}, &j\equiv 1\pmod 2,\\
\end{cases}
$$
where indices on the right-hand side are reduced modulo $6$.
A slightly tedious calculation now shows that $\varphi$ is indeed an adjacency preserving bijection of the vertex sets. In similar fashion, we could show the other two graphs above are isomorphic in case when $n/2$ is odd.

\medskip
A slightly different situation happens for $X\cong\HTG(m,n,\ell)$, when $n=6$ and $m\ne 3$. In this case, each vertical edge is contained in $3$ different $6$-cycles, while each horizontal edge and each jump is contained in $2$ different $6$-cycles. In this case, we have $\signat(X)=(2,2,3)$ and $X\cong\HTG(m,6,0)$ for $m$ even, or $\HTG(m,6,3)$ for $m$ odd. Again we can also prove the isomorphism relations by defining an apropriate bijection, but we omit the details.
\medskip

The remaining graphs of girth $6$ now have no other $6$-cycles but those with 4 vertical and 2 horizontal edges or jumps, implying that 
their signature is $\signat(X)=(2,2,2)$. Among these graphs, we note that for all $n\equiv 0\pmod 4$, graph $\HTG(2,n,\frac{n}{2})$ is isomorphic to the generalized Petersen graph $\GP(n,\frac{n}{2}-1)$ of girth $6$. The appropriate graph isomorphism identifies the vertices in the first column $v_{0,j}$ of $\HTG$ with vertices $u_j$ on the outer cycle of the generalized Petersen graph, and reorders the vertices of the second column by mapping $v_{1,j}$ to $v_{j(\frac{n}{2}-1)}$, hence mapping edges between columns in HTG into spokes of GP graph, as can be checked by some effort. The remaining 3-RDR HTG graphs are not isomorphic to any GP graphs, as their signatures are different.
\end{proof}
\subsection{\sc A new family of 3-RDR Cayley graphs}\label{subs Xn}

As first observed by computations in Magma, all vertex-transitive 3-RDR graphs of orders up to 30 are isomorphic to generalized Petersen graphs or honeycomb toroidal graphs, see Table~\ref{tab: ord<36}. 
However, there are two  graphs not belonging to any of these two families among the graphs of order $36$.
One of them belongs to a new family of vertex-transitive 3-RDR graphs that we will now construct as undirected Cayley graphs over group $G=S_3\times D_n$, where $D_n=\langle r,z\mid r^n=z^2=(zr)^2=1\rangle$ is the dihedral group of order $2n$, $n\geq 3$.
Let $a=((12),1)$, $b=((13),z)$ and $c=((23),zr)$ be three involutions from group $G$, and let $S=\{a,b,c\}$.
Now define graphs $X_n=\Cay(G,S)$ with vertex set $V(X_n)=G$ and edge set  $E(X_n)=\{\{g,gs\}: g \in G,s\in S\}$.
\begin{theorem}\label{thm Xn}
Graph $X_{n}$, $n\geq 3$, is a connected vertex-transitive 3-RDR graph of order $12n$, and the following holds:
\begin{enumerate}[(i)]
\item $\girth(X_3)=6$ and $\signat(X_3)=(0,1,1)$.
\item $\girth(X_n)=8$, $\signat(X_n)=(5,5,6)$ for $n\geq 4$, with
$$X_n\cong 
\begin{cases}
\GP(6n,2n-1),& n\equiv 1\pmod 3,\\
\GP(6n,2n+1),& n\equiv -1\pmod 3,\\
\end{cases}$$
and $X_n$ is not isomorphic to any generalized Petersen graph for $n\equiv 0\pmod 3$.
\end{enumerate}
\end{theorem}
\begin{proof}
Observe first that the connection set $S$ contains three involutions. Since $bab=((23),1)$ and $cac=((13),1)$, we easily see that $S$ generates G, hence $\Cay(G,S)$ is connected and cubic of order $|G|=12n$.
Now let $H=A_3\times D_n$ be a subgroup of index $2$ in $G$, consisting of elements $h=(h_1,h_2)$ with first entry $h_1$ an even permutation. Then
for any $h\in H$, elements $ha,hb,hc$ have an odd permutation as the first entry, hence $ha,hb,hc\notin H$. Therefore, cosets $H$ and $Ha$ are bipartition sets of vertices of graph $X_{n}$ and a natural 3RD-coloring $f$ can be defined by 
$$f(h_1,h_2)=\begin{cases} 
\emptyset, &h_1\in A_3,\\
\{1\}, &h_1=(12),\\
\{2\}, &h_1=(23),\\
\{3\}, &h_1=(13).
\end{cases}$$
This shows that $X_{n}$ is a vertex-transitive 3-RDR graph.

Since $X_n$ is bipartite, its cycles are of even length. We compute the 24 reduced $S$-words of length $4$ in $G$ to obtain 16 distinct elements as follows: 
\begin{align*}
abab=caca, 
abac=baca, &
abcb=bcba, 
acab=caba,\\
acac=baba,
acbc=cbca, &
babc=cbab,
bcab=cacb,\\
abca, acba, bcbc, bcab, & bacb, cbac, cbcb, cabc.
\end{align*}
As none of these elements is equal to $1_G$, there are no $4$-cycles in the graph.
Since $bc=((321),r)$, we have that $(bc)^3=1$ iff $n=3$. In this case, equality $bcbcbc=1$ represents a $6$-cycle in $X_3$, thus $\girth(X_3)=6$. Moreover, there is a unique $6$-cycle containing edge $e_b=\{1_G,b\}$ and edge $e_c=\{1_G,c\}$, and there are no $6$-cycles containing edge $e_a=\{1_G,a\}$, so $\signat(X_3)=(0,1,1)$. 

For $n\geq 4$, reduced words of length $4$ are all distinct from reduced words of length $2$ $ab,ac,ba,bc,ca,cb$, hence there are no $6$-cycles either. However, $8$-cycles exist since equality $abab=caca$ yields  $ababacac=1$, etc.
Again, observe the edge $e_a\in E(X_{n})$. Equations above imply
$$ababacac=abacacab=abcbabcb=acababac=acacabab=acbcacbc=1,$$ 
so edge $e_a$ is contained in exactly $6$ different $8$-cycles starting at vertex $1_G$.
 In similar fashion, we see that edges $e_b=\{1_G,b\}$ and $e_c=\{1_G,c\}$ are contained each in exactly $5$ different $8$-cycles starting at vertex $1_G$. Thus we have $\signat(X_n)=(5,6,6)$ for $n\geq 4$.
 
For the isomorphisms, suppose first that $n\equiv -1\pmod 3$. To show that $X_n\cong \GP(6n,2n+1)$, observe that element $bc=((321),r)\in G$ has order $3n$, so equality $bcbc\cdots bc=1$ gives an outer cycle of length $6n$ in $X_n$. We label the vertices of $G$ by $u_i,v_i$, $i\in{\mathbb Z}_{3n}$ as
$$\begin{array}{llllll}
u_0=1,& u_1=b,& u_2=bc,& u_3=bcb,&\ldots ,&u_{3n-1}=bc\cdots b,\\
v_0=a,& v_1=ba,& v_2=bca,& v_3=bcba,&\ldots ,&v_{3n-1}=bc\cdots ba.
\end{array}$$
Obviously, the cycle $1,b,bc,\ldots, bc\cdots b \in X_n$ corresponds to the outer cycle of $\GP(6n,k)$ and edges $\{g,ga\}$ in $X_n$ correspond to spokes $\{u_i,v_i\}$. Now we show that any edge $\{g,gb\}$, where $g=v_i$ is in the bottom row, corresponds to edge $\{v_i,v_{i+k}\}$ for $k=2n+1$.
Without loss of generality, it is enough to check this for $i=0$. Let $g=v_0=a=((12),1)$. Then $gb=ab=((321),z)$, and since $n\equiv -1\pmod 3$, we also have  $$v_{k}=v_{2n+1}=b(cb)^n a=((13)(123)^n(12),zr^{-n})=((321),z).$$
Hence edge $\{a,ab\}$ in $X_3$ corresponds to edge $\{v_0,v_{2n+1}\}$ in $\GP(6n,2n+1)$.
In similar fashion, one can show that $X_n\cong \GP(6n,2n-1)$ when $n\equiv 1\pmod 3$.
Finally, for $n\equiv 0\pmod 3$, $n\geq 4$, vertex-transitive $3$-RDR generalized Petersen graphs of appropriate order and girth have signature $(2,2,4)$ by Theorem \ref{3RDR GP}, so $X_n$ is not isomorphic to any of them.
\end{proof}
Since HTG graphs have girth $\leq 6$ and signatures as in Theorem~\ref{thm HTG}, we also have:

\begin{corollary}\label{Xn ni HTG}
For $m\geq 1$, graphs $X_{3m}$ are an infinite family of vertex-transitive 3-RDR graphs of order $36m$, which are not isomorphic to $\HTG$ or $\GP$ graphs.
\end{corollary}

We note that among the vertex-transitive  3-RDR graphs of orders up to $36$ (Table \ref{tab: ord<36}), there is a unique graph which is not a $\HTG$, $\GP$ or $X_n$ graph. By classification of vertex-transitive cubic graphs of girth 4 by \v Sparl, Eiben and Jajcay \cite{Sparl Eiben Jajcay}, all such graphs must appear as prisms, generalized prisms, M\"obius ladders or generalized graph truncations of arc-transitive graphs of degree $4$ by the $4$-cycle. There are several more unidentified VT 3-RDR graphs of girths $4,6,8,10$ with orders $\geq 48$, see last column of Table~\ref{tab: bicubic}.

\section{Concluding remarks}\label{sec: remarks}
As our data shows, the project of classifying all vertex-transitive $3$-RDR graphs is far from finished. Identifying some other infinite families of vertex-transitive $3$-RDRs as Cayley graphs over appropriate groups, or alternatively, as covering graphs over $K_{3,3}$ seems possible. Using the recent results on graphs on vertex-transitive graphs of small girths  from \cite{Sparl Eiben Jajcay} and \cite{Potocnik Vidali girth six}, at least the graphs of girth $\leq 6$ might be possible to classify completely. However, we know that unidentified vertex-transitive $3$-RDRs of girth $8$, $10$ and possibly even larger exist for orders $60,72,84,96,\ldots$. 

On the other hand, our current investigations of $d$-RDR graphs also leave several other questions open. For instance, all $3$-RDR graphs of small orders that satisfy criteria from Theorem \ref{krit1} also satisfy criteria from Theorem \ref{krit2}. Could this be true in general, and for all vertex-transitive $d$-RDR graphs? We were also not able to identify any $3$-RDR graphs not satisfying criteria from Theorem \ref{krit2}, although they might exist. We also note that all vertex-transitive $3$-RDR graphs identified in our investigations admit a regular group of automorphisms and hence are Cayley graphs. Are there also any non-Cayley vertex-transitive $3$-RDRs (or $d$-RDRs)? Other symmetry properties could also be investigated in the context of $d$-RDR graphs. For instance, for graphs of orders up to $96$, the Tutte-Coxeter graph of order 30 is the only arc-transitive bicubic graph that is not $3$-RDR. We look forward to answering any of these questions in further investigations.

\section{Acknowledgments}
This work is supported in part by the Slovenian Research Agency ARIS, research program P1-0285 and research projects J1-3001, J1-1694, J1-2451. 

\end{document}